\documentclass[preprint,3p,authoryear,times]{elsarticle}

\journal{} 
\usepackage[fleqn]{amsmath}
\usepackage[]{amsfonts}
\usepackage[]{amsthm}
\usepackage[]{amssymb}
\usepackage{empheq}
\usepackage{lscape} 

\usepackage[mathlines]{lineno}

\usepackage[]{algorithm2e}
\usepackage[english]{babel}
\usepackage[T1]{fontenc}

\usepackage{hyperref} 
\usepackage{graphicx}
\usepackage{caption,subcaption}
\usepackage{tikz}
\usepackage{pstricks}
\usepackage{xstring}
\usepackage{multirow}
\usepackage{hhline}

\definecolor{darkblue}{rgb}{0,0,0.5}

\hypersetup{colorlinks=true,linkcolor=black,citecolor=darkblue,urlcolor=darkblue} 

\newcommand{\newref}[2]{\hyperref[#2]{#1~\ref*{#2}}} 

\usepackage{longtable}
\usepackage{pdflscape}
\usepackage{tabularx}
\usepackage[toc,page]{appendix}

\SetArgSty{textnormal}
\usepackage{float}
\usepackage{booktabs}
\usepackage[final]{changes}
\definechangesauthor[color=orange]{rev1}
\definechangesauthor[color=red]{rev2}

\usepackage{rotating} 
\usepackage{diagbox} 

\begin{document}
\begin{frontmatter}

\title{Multimodal urban transportation network equilibrium including intermodality and shared mobility services}
\author[1]{Khadidja Kadem\corref{cor}}
\cortext[cor]{Corresponding author.} 
\ead{khadidja.kadem@univ-eiffel.fr}
\author[1,2]{Mostafa Ameli}
\author[1]{Mahdi Zargayouna} 
\author[1]{Latifa Oukhellou}

\address[1]{COSYS-GRETTIA, Univ Gustave Eiffel, Marne-la-Vallée, F-77454, France}
\address[2]{Department of Electrical Engineering and Computer Sciences,  University of California, Berkeley, Berkeley, USA} 

\begin{sloppypar} 
\begin{abstract} 

Shared Mobility Services (SMSs) are transforming urban transportation systems by offering flexible travel options. These services, which help reduce the number of cars on the roads, have the potential to enhance the transportation system's performance, leading to improvements in travel times and emissions. This emphasizes the importance of assessing their impact on the system and users’ choices, particularly when integrated into complex multi-modal systems that include public transport (PT). However, many studies overlook the synergies between SMSs and PT, leading to inaccurate traffic estimations and planning. This research presents an extensive review of multi-modal transportation system models incorporating SMSs. It then introduces a multimodal traffic assignment model including almost all mobility options in urban transportation systems applicable in both continuous and integer settings, leading to a Mixed-Integer Bilinear Programming (MIBLP) formulation. The model comprises diverse travel options, including SMSs, and accounts for intermodality by allowing commuters to combine modes to optimize time and monetary expense. An in-depth examination of commuters' mode and path choices on two test cases and an analysis of the price of anarchy reveals the disparities between user equilibrium and system optimum in such complex networks.

\end{abstract}

\begin{keyword}
ridesharing \sep carpooling \sep public transport \sep traffic assignment \sep system optimum \sep user equilibrium \sep Intermodality

\end{keyword}

\end{sloppypar} 
\end{frontmatter}


\begin{sloppypar} 

\section{Introduction}

With the rapid growth of demand for urban transportation systems, addressing the population's travel needs has become challenging~\citep{nalin2025evaluation}. Because expanding the physical transportation network is not always feasible and beneficial, new mobility services have emerged, transforming and revolutionizing how people move within cities~\citep{du_modeling_2022}. In particular, new shared mobility services (SMSs) have gained popularity due to their cost-effectiveness and potential to extend transit coverage in underserved areas~\citep{whitmore2022integrating}. By promoting shared mobility and reducing the number of cars in the network, these services can potentially reduce traffic congestion and emissions and alleviate parking space demands~\citep{SCHALLER20211}. Furthermore, on-demand SMSs can radically improve mobility, especially in low-demand and low-density regions where public transport (PT) is unable to provide first and last mile transportation (e.g., from origin to the transportation hub or from the hub to destination)~\citep{boarnet_firstlast_2017}. These on-demand services are extremely convenient for commuters, enabling them to request a ride precisely when and where they need it. The e-hailing service is a prime example of redefining traditional taxi services. This convenience also offers new opportunities to enhance the transportation system efficiency by exploiting the sharing of the same vehicle by multiple passengers, giving rise to the ridesharing service. In the same context, carpooling emerges as another shared mobility service, where commuters with overlapping itineraries and similar schedules share a privately owned car for the trip~\citep{murphy_shared_2016}. The main difference between ridesharing and carpooling is the fact that in carpooling, the driver is also a commuter, with specific travel needs. They accept to share their car with passengers for a part of their trip in order to reduce their travel costs. In ridesharing, and similar to a shared taxi service, the driver is performing the service to a pool of passengers as a job. Moreover, e-hailing and ridesharing, though operated by the same provider, follow distinct mechanisms. Ridesharing introduces greater complexity due to the shared nature of trips, impacting the matching process, costs, waiting times, as well as vehicle management since it involves handling partially occupied vehicles.

Integrating these new SMSs with PT has shifted traveling behavior and created a more complex multi-modal transportation system~\citep{pi_general_2019}. On the one hand, commuters face mode and path choices, aiming to reduce travel expenses. On the other hand, the system owners or traffic engineers strive to optimize mobility services and meet demand while minimizing overall costs \citep{koohathongsumrit2024multimodal}. This needs a complete understanding of the various transportation services accessible and their associated costs.

Numerous studies have focused on modeling multi-modal transportation systems. Recently, many works have incorporated SMSs into their models. However, a less explored aspect concerns modeling and evaluating combinations of multiple mobility options to accomplish one trip. This concept, known as intermodality, assumes that a commuter chooses a set of modes to accomplish their trip and the path and location to transfer between the modes. Particularly, using intermodality involving public transport and SMSs for first and last mile transportation is seldom explored. This lack of integration results in inaccurate traffic predictions and planning. Consequently, individuals may not be able to fully optimize their travel choices, leading to potential inefficiencies~\citep{zhu_analysis_2020}.

Addressing this gap requires the development of comprehensive transportation models that consider all travel options available to commuters. By integrating SMSs and PT into a unified framework and providing accurate information on costs and travel times, the travel and service expenses of both commuters and the urban transportation system can be optimized. Two traffic assignment principles are used to formulate the commuter's decisions: User Equilibrium (UE) and System Optimum (SO). In UE, known as the first Wardrop principle~\citep{wardrop_theoretical_1952}, every commuter tries to minimize their own cost, and equilibrium is reached when no one is willing to change their choices to achieve lower costs. The second Wardrop principle defines SO as a system-oriented equilibrium and assumes that the system's total cost is minimized~\citep{wardrop_theoretical_1952}. Analyzing the gap between these two principles, referred to as the Price of Anarchy (PoA), will allow the measurement of the system's degradation due to its agents' selfish and conflicting behavior. 

\subsection{Literature Review}

Several studies in the literature addressed the modeling of shared mobility services, including background traffic~\citep{zhu_analysis_2020}. \cite{xu_complementarity_2015} and~\cite{di_ridesharing_2017} studied the impacts of the carpooling service on traffic congestion. \cite{di_unified_2019} proposed a UE framework to analyze the congestion effects of carpooling and e-hailing services. \cite{ma_ridesharing_2020} formulated a UE model with the carpooling service and studied the effects of an OD-based surge pricing strategy. \cite{li_path-based_2020}, \cite{zhong_dynamic_2020} and~\cite{wang_convex_2021} also considered the carpooling service in a UE framework. Similarly,~\cite{yan_stochastic_2019} studied carpooling demand, with the same origin and destination for riders and drivers, under the Stochastic User Equilibrium (SUE) principle, in which commuters are assumed to minimize their perceived travel costs. They showed that carpooling is preferred for long trips. For the e-hailing and ridesharing services,~\cite{nourinejad_ride-sourcing_2020} considered the e-hailing service while~\cite{beojone_inefficiency_2021} analyzed the efficiency of both services in reducing traffic congestion. They showed that fleet size is an important factor for ridesharing performance due to the congestion effects of empty fleet vehicles. In the same context,~\cite{alisoltani_can_2021} showed that the ridesharing service can reduce traffic congestion in large-scale urban transportation systems. Still, its impacts are not significant in small and medium-scale cities. However, these studies did not consider public transport and its implications on passengers' behavior and the state of road traffic.

Some works have included both SMSs and PT in their modeling. For example,~\cite{qian_modeling_2011} studied the use of transport modes and network performance when the carpooling service is integrated with public transport. Similar work was done by~\cite{sun_multi-class_2021}, in which the authors formulated a stochastic user equilibrium (SUE) with multiple classes of user demand. In \cite{yao_general_2023}, the authors formulated a UE model integrating the matching decisions of the carpooling platform with the mode choices of commuters among solo-driving, carpooling and metro. \cite{zhang_pool_2021} analyzed the e-hailing and ridesharing pricing strategies in a multi-modal network, considering PT. \cite{liang_dynamic_2023} proposed a UE framework, under time uncertainties, to analyze the commuters' mode and departure time choices. \cite{narayan_does_2019} studied the impacts of the e-hailing service on the PT demand and use of private cars in Amsterdam, Netherlands. They showed that e-hailing could replace private cars and attract PT demand with a larger fleet. \cite{wei_modeling_2020} modeled the carpooling service in a multi-modal single-region network and proposed two pricing schemes to alleviate traffic congestion. \cite{tang_multi-modal_2021} investigated the integration of shared autonomous vehicles into the multi-modal network, considering PT and parking space constraints. These studies only considered railway transit, disregarding interactions between the road network and public transport. To this end,~\cite{fayed_utilization_2023} investigated the use of underutilized bus lanes for the ridesharing service and its impacts on the solo-pooled rides split. However, in these studies, the commuters make single-mode trips and cannot switch modes between their origin and destination, i.e., intermodality is not considered. This results in the use of PT being under-estimated or over-estimated, as in most cases, PT does not cover the first and last mile of the journey. Thus, intermodality can be divided into two categories in the literature. In the first one, commuters use personal means of transport (i.e., car, bike, or walking) for the first or last mile of the trip, while in the second category, SMSs are used for this purpose as a complement to public transport. 

Various studies have considered the first category of intermodality by exclusively targeting the use of cars to reach a PT hub, known as Park-and-Ride mode in the literature \citep{ye_joint_2021, zheng_multimodal_2020, liu_doubly_2017, fan_bilevel_2014}. The use of bikes to complete the first mile and last mile has also been intensively investigated (See, e.g.,~\cite{wang_designing_2022, geurs_multi-modal_2016}). However, considering intermodality with SMSs adds another level of complexity. Allowing SMSs to be used as a complementary service instead of an exclusive door-to-door service requires the adaptation of the model to ensure that the supply and demand flows remain coherent. In addition, it allows more opportunities for matching since the demand, departing (or exiting) at transfer points also needs to be considered. Thus, only a few recent studies have targeted intermodality with SMSs and PT. \cite{ke_equilibrium_2021} analyzed the competitiveness and complementarity of e-hailing with PT. They showed that fleet size and trip fares impact the use of e-hailing services as a first-mile and/or last-mile option. \cite{zhang_integrating_2021} and~\cite{liu_integrating_2023} formulated the SUE when the e-hailing service is complementary to PT.
Similarly,~\cite{zhu_analysis_2020} investigated the impacts of e-hailing and ridesharing services on PT ridership and showed that ridesharing can be a competitor of PT when low fares are applied. In~\cite{du_modeling_2022}, the authors proposed an SUE multi-modal framework where commuters access metro stations via bus, car, or e-hailing. The mode and path choices are handled with a cross-nested logit model. The carpooling mode was also considered for both drivers and passengers, exclusively as a door-to-door service. In~\cite{pi_general_2019}, commuters can choose bi-modal travel options with solo driving (park-and-ride mode), carpooling, and e-hailing to access transfer stations. Travel behavior is modeled via a multi-layer nested logit model, while the matching for the carpooling service is not considered with different OD pairs between the passenger and the driver.

\newref{Table}{tab:1} presents a characterization of the existing models in the literature. The traffic assignment column categorizes the studies based on equilibrium principles and highlights the lack of work investigating the difference between UE and SO in such a complex multi-modal system. More specifically, existing works in the literature targeting the Price of Anarchy (PoA) concept only focus on a single mode. Most often, they provide an experimental and theoretical analysis of PoA for the car mode (e.g.,~\cite{colini-baldeschi_when_2020} and~\cite{qiao_price_2023}). Other studies, such as~\cite{fielbaum_how_2022} and~\cite{chau_quantifying_2018}, examine the PoA in a system with one mobility service in which passengers share a vehicle. However, these works did not address PoA in a transportation system with the presence of public transport and shared mobility services.

The generalized cost column presents the terms included in the travel costs: travel time stands for in-vehicle time, waiting time for the time spent before the vehicle arrives at the passenger's position, the service time is the time to be matched, to get in and out of the vehicle and the time spent at PT stops, the monetary cost is the cost of fuel and trip fares. In this work, we choose to focus on these four costs as the major costs impacting the commuters' behavior. However, the existing works may consider other costs, such as inconvenience or detour costs for SMSs (e.g. \cite{di_unified_2019}, \cite{li_path-based_2020}, \cite{wang_convex_2021}, and \cite{yao_general_2023}). These costs can be implicitly incorporated in the service time as defined in this paper. 
Lastly, the travel options column denotes the categories of travel modes and their travel options considered in each model. There are four categories: Personal means, Public Transport (PT), Shared Mobility Services (SMSs), and Intermodality, which include park-and-ride, walking and/or biking combined with PT (W/B), and using SMSs with PT. To the best of our knowledge, only the works of~\cite{du_modeling_2022} and~\cite{pi_general_2019} present models, considering a single type of SMS and intermodality. However, the ridesharing service, which requires additional constraints compared to the carpooling service, was not considered. These constraints concern the conservation of the fleet size, the availability of cars to perform the service, and the states conversion (idle to occupied and vice-versa). Moreover, matching passengers with different OD pairs, and with different service order (FIFO vs LIFO) needs to be addressed in ridesharing. None of the mentioned works of the literature have provided a complete approach to address all aspects simultaneously.

\begin{table}[!]
\centering
\caption{\textbf{Comparison of the different studies on modeling of shared mobility services}}
\resizebox{1\columnwidth}{!}{
\begin{tabular}{c|ccc|cccc|ccc|cc|ccc|ccc}

\multirow{4}{*}{Research}  & \multicolumn{3}{c|}{{\begin{tabular}[|c|]{@{}c@{}} Traffic \\ Assignment \end{tabular}}} & \multicolumn{4}{c|}{{\begin{tabular}[|c|]{@{}c@{}} Generalized \\ Cost \end{tabular}}} & \multicolumn{11}{c}{Travel Options}  \\
\cline{9-19} 
& & & & & & & & \multicolumn{3}{c|}{Personal Means} & \multicolumn{2}{c|}{PT} & \multicolumn{3}{c|}{SMSs}  & \multicolumn{3}{c}{Intermodality} \\
\cline{2-19} 
& \begin{tabular}[c]{@{}c@{}}User\\ Equilibrium \end{tabular} & \begin{tabular}[c]{@{}c@{}} System \\ Optimum \end{tabular}  & \begin{tabular}[c]{@{}c@{}} Stochastic \\ UE \end{tabular}  & \begin{tabular}[c]{@{}c@{}} Travel \\ Time \end{tabular}  & \begin{tabular}[c]{@{}c@{}} Waiting \\ Time \end{tabular}  & \begin{tabular}[c]{@{}c@{}} Service \\ Time \end{tabular}  & \begin{tabular}[c]{@{}c@{}} Monetary  \\ Cost \end{tabular} & Car & Walking & Biking & Bus & Metro & Carpooling & E-hailing & \begin{tabular}[c]{@{}c@{}} Ride \\ Sharing \end{tabular}  & \begin{tabular}[c]{@{}c@{}} Park \& \\ Ride \end{tabular}  & \begin{tabular}[c]{@{}c@{}} W/B \\  with PT \end{tabular} & \begin{tabular}[c]{@{}c@{}} SMSs  \\ with PT \end{tabular}\\ 
\hline
\cite{qian_modeling_2011}&x& & &x& &x&x&x& & & &x&x*& & & & & \\
\cite{fan_bilevel_2014}& & &x&x&x&x&x&x& & & &x& & & &x& & \\
\cite{xu_complementarity_2015} & x & & & x & x & & x& x& & & & & x& & & & & \\
\cite{geurs_multi-modal_2016}& & & &x&x& &x&x&x&x&x&x& & & &x&x& \\
\cite{di_ridesharing_2017} & x & & & x & x & & x& x& & & & & x& & & & & \\
\cite{liu_doubly_2017}&x& & &x&x& &x&x& & & &x& & & &x& & \\
\cite{di_unified_2019} & x & & & x & x & & x & x & & & & & x & x & & & & \\
\cite{narayan_does_2019}&x& & &x&x& &x&x&x&x& &x& &x& & & & \\
\cite{pi_general_2019}&x& & &x&x& &x&x& & &x&x&x*&x& &x& &x \\
\cite{yan_stochastic_2019} & & &x&x&x& &x&x& & & & &x*& & & & & \\
\cite{li_path-based_2020}&x& & &x& & &x&x& & & & &x& & & & & \\
\cite{ma_ridesharing_2020}&x& & &x& & &x&x& & & & &x& & & & & \\
\cite{nourinejad_ride-sourcing_2020}& & & &x&x& &x& & & & & & &x& & & & \\
\cite{wei_modeling_2020}& & & &x&x& &x&x& & & &x&x*& & & & & \\
\cite{zheng_multimodal_2020}& & &x&x&x& &x&x& & & &x& & & &x& & \\
\cite{zhong_dynamic_2020}&x& & &x&x& &x&x& & & & &x& & & & & \\
\cite{zhu_analysis_2020}&x& & &x&x& &x& & & & &x& &x&x*& & &x \\
\cite{alisoltani_can_2021}& & & &x&x&x& &x& & & & & &x&x& & & \\
\cite{beojone_inefficiency_2021} & & & &x&x& & & x & & & & & &x&x& & & \\
\cite{ke_equilibrium_2021}& & & &x&x& &x& &x&x& &x& &x& & &x&x \\
\cite{sun_multi-class_2021}& & &x&x&x& &x&x& & & &x&x& & & & & \\
\cite{tang_multi-modal_2021}& &x& &x& &x&x&x& & & &x& & &x*& & & \\
\cite{wang_convex_2021}&x& & &x& &x&x&x& & & & &x& & & & & \\
\cite{ye_joint_2021}& & &x&x& & & &x&x&x& &x& & & &x& & \\
\cite{zhang_integrating_2021}& & &x&x&x& &x& & & & &x& & & & &x&x\\
\cite{zhang_pool_2021}& & & &x&x& &x& & & & &x& &x&x& & & \\
\cite{du_modeling_2022}& & &x&x& & &x&x& & &x&x&x&x& &x& & \\
\cite{wang_designing_2022}& & & &x&x&x&x&x&x&x&x& & &x& & &x& \\
\cite{fayed_utilization_2023}& & & &x&x&x&x&x& & &x& & &x&x& & & \\
\cite{liang_dynamic_2023}&x& & &x& & &x&x& & & &x& &x&x*& & & \\
\cite{liu_integrating_2023}& & &x&x&x& &x&x& & & &x& &x& & & &x\\
\cite{yao_general_2023}&x& & &x& & &x&x& & & &x&x& & & & &  \\
\textbf{This study} &\textbf{x}&\textbf{x}& &\textbf{x}&\textbf{x}&\textbf{x}&\textbf{x}&\textbf{x}&\textbf{x}&\textbf{x}&\textbf{x}&\textbf{x}&\textbf{x}&\textbf{x}&\textbf{x}&\textbf{x}&\textbf{x}&\textbf{x} \\
\hline 
\multicolumn{19}{l}{\small{*: The sharing is for the same Origin-Destination (OD) pair}}
\end{tabular} \label{tab:1}
}
\end{table}

\subsection{Contribution}

One of the goals of this study is to fill the research gap highlighted in \newref{Table}{tab:1} by proposing a comprehensive model for multi-modal transportation systems, integrating SMSs with PT, and enabling intermodality. We propose a combined framework for carpooling, e-hailing, ridesharing. We aim to address UE and SO principles in a tractable (analytical) framework, thus eliminating the need to calibrate choice models for commuters. We address the equilibrium problem with the Beckmann transformation and prove that the solution of our model satisfies the UE conditions in the continuous setting, with respect to a delay induced by the use of SMSs, and subject to capacity constraints. Therefore, our general link-based formulation of the mode and path choice under the UE principle can be converted to other types of equilibria to represent the commuters' choices. Moreover, we represent traffic flow using integer variables and examine the theoretical implications of this approach on the UE conditions of such multimodal programs. This choice is mainly motivated by the discrete nature of traffic flows, especially in SMSs contexts, where individual or vehicle units are involved. Representing traffic as integers allows for more accurate modeling of these dynamics, including intermodal interactions and the passenger-driver matching process, thereby aligning the analysis more closely with real-world scenarios. However, our general formulation of the objective function and constraints ensures that our model remains valid with continuous flow variables.

The proposed model enables the investigation of different planning scenarios by providing network traffic distribution. 
In this study, we perform a comparative analysis of the two equilibria to investigate the impact of intermodality on the transportation system performance. In particular, we compare the system, including PT and SMSs, competing with the cooperative scenario wherein intermodality is allowed.

Hence, the contributions of this paper can be summarized as follows: 
\begin{itemize}
  \item It establishes a comprehensive modeling approach for multi-modal transportation networks considering personal transportation means (car, walking, and biking), public transport (bus and metro), as well as shared mobility services (carpooling, ridesharing, and e-hailing).
  \item It explicitly formulates interactions between shared mobility services and public transport regarding congestion effects and intermodality. It also considers the park-and-ride mode and the use of walking and biking for first-mile and last-mile transportation.
  \item It calculates the user equilibrium via the Beckmann formulation and the system optimum when the model simultaneously handles mode and path choices and passenger-driver matching. This model is applicable in both continuous and integer settings, resulting in a Mixed-Integer Bilinear Programming (MIBLP) formulation.
  \item The proposed model provides insights into factors impacting the use of the modes in a synthetic multi-modal network. It also compares user equilibrium and system optimum solutions by analyzing the price of anarchy in such multi-modal systems.
  \item The Sioux Falls urban transportation is also considered to analyze commuters’ behavior and scale up the proposed model's applicability to realistic scenarios.
\end{itemize}

The rest of the paper is organized as follows. \newref{Section}{sec::model} presents the mathematical formulation for the proposed framework. In \newref{Section}{sec::ER}, we first analytically validate our formulation on a synthetic network. Then, we extend the application to a larger network to provide more insights on mode sharing and intermodality of a complete multi-modal urban network. We explore different transportation scenarios and provide an in-depth analysis of each. Additionally, we discuss the computational performance of our model and its scalability to large-scale networks. Finally, we present the concluding remarks in \newref{Section}{sec::conclusion}.

\section{Model formulation}
\label{sec::model}
This section presents the modeling approach for multi-modal urban transportation networks considering shared mobility services integrated with public transport. For clarity's sake, we first present a simplified version of the model without intermodality. Then, we reformulate the model's constraints to consider intermodality. The list of notations used in this paper is presented in \newref{Table}{tab:2}.

\begin{table}[!]\centering
\caption{\textbf{List of notations.}}
{\footnotesize
\begin{tabular}{p{1.2cm}p{14cm}}\hline
\multicolumn{2}{l}{\textbf{\underline{Network Structure}} } \\
E & Set of nodes. \\
O & Set of origin nodes. \\
D & Set of destination nodes. \\
TR & Set of transfer nodes. \\
N & Set of networks ; $N = \{Road (RN),\; Metro (MN),\; Bike (BN),\; Walk (WN)\}$. \\
$A_n$ & Set of links in network \textit{n.}\\
$A$ & Set of links: $A = \cup_{n \in N} A_n$.\\
$\Psi$ & Set of all travel modes ; $\Psi$ = \{ car, bus, M, CP, CD,
EH, RS, W, B, $I_{m1,m2} $ \} 
where $I_{m1,m2}$ is a combination of mode $m1$ and $m2$. \\
V & Set of fleet vehicles for ridesharing and e-hailing services.\\
$P_{ij}^m$ & Set of paths for origin-destination (OD) pair $(i,j)$ with mode $m$. \\
$o_a$ & Binary coefficient equal to 1 if link \textit{a} is an origin link ; 0 Otherwise. This parameter depends on the network structure and the OD demand. \\
$d_a$ & Binary coefficient equal to 1 if link \textit{a} is a destination link ; 0 Otherwise. This parameter depends on the network structure and the OD demand. \\
$\delta_{a,p,m}^{(i,j)}$ & Binary coefficient equal to 1 if path $p$ of mode $m$ traverse link \textit{a} between OD pair $(i,j)$ ; 0 Otherwise. \\
$q^{(i,j)}$ & Travel demand for OD pair $(i,j)$. \\
$Q$ & Total demand for the network.
\\ 
\textbf{\underline{Indices} }& \\
$i,j,r,s$ & Index of node, $i,j,r,s \in E$. \\
$a$ & Index of link (edge), $a \in A$. \\
$m$ & Index of mode, $m \in \Psi$. \\
$p$ & Index of path, $p \in P_{ij}^m $. \\
$n$ & Index of network,  $n \in N$. \\
\\
\multicolumn{2}{l}{\textbf{\underline{Input Parameters  - [Value for the experiments]}} } \\
$L_a$ & Length of link \textit{a}. \\
$Sp_n$ & Mean speed in network $n$  - [M: 60, B: 10, W: 3]. \\
$freq_{m, a}$ & Frequency of travel mode $m$ on link \textit{a} (vehicles/time unit)  - [Bus: 3, M: 6]. \\
$R_m$ & Meeting rate of travel mode $m$ (service/time unit)  - [CP: \{synthetic network: 100 , Sioux Falls: 500\}, RS: \{synthetic network: 200 , Sioux Falls: 1000\}, EH: \{synthetic network: 200 , Sioux Falls: 1000\}]. \\
$S_m$ & Average service time for mode $m$  - [Bus: 0.04, M: 0.02, CP: 0.04, RS: 0.05, EH: 0.03]. \\
$P_m$ & Average parking time for mode $m$  - [car: 0.17, B: 0.08]. \\
$TF_{m,a}$ & Trip fare of travel mode $m$ on link \textit{a} - [Bus: 0.3, M: 0.3, CP: 0.7, CD: 0.7, RS: 0.9, EH: 1.1]. \\
$PF_{m}$ & Parking fare for travel mode $m$  - [car: 1]. \\
$CAP_m$ & Maximum passenger capacity for the travel mode $m$  - [CP: 1 , RS: 2]. \\
$Pk\_cap_s$ & Parking capacity for transfer node $s$  - [synthetic network: 100, Sioux Falls: 800]. \\
$|V|$ & Fleet size for the ridesharing and e-hailing services  - [synthetic network: 1000, Sioux Falls: 3500]. \\
$\alpha$ & Monetary cost per unit of time (monetary unit/time unit)  - [synthetic network: 5, Sioux Falls: 20]. \\
$\gamma$ & Monetary cost per unit of distance (monetary unit/distance unit)  - [0.25].\\
\\
\textbf{\underline{Variables}} &  \\
$t_{a,n}$ & Travel time on link \textit{a} belonging to network \textit{n}. \\
$WT_{m,a}$ & Waiting time for travel mode $m$ on link \textit{a} . \\
$ST_{m,a}$ & Service time for travel mode $m$ on link \textit{a}. \\
$C_{m,a}$  & Monetary cost of using travel mode $m$ on link \textit{a}. \\
$x_a$ & Aggregated traffic flow on link \textit{a}, (Integer variable). \\
$x_{a,m}$ & Aggregated traffic flow on link \textit{a} with travel mode $m$, (Integer variable).\\
$q_m$ & Travel demand for mode $m$, (Integer variable). \\
$q_m^{(i,j)}$ & Travel demand for mode $m$ and OD pair $(i,j)$, (Integer variable). \\
$q_t^{(i,j)}$ & Number of empty or occupied service vehicles between OD pair $(i,j)$ ; $t \in \{e,o\}$, (Integer variable). \\
$q_o^{(i,j)\;r}$ & Number of occupied service vehicles between OD pair $(i,j)$ stopping at node $r$, (Integer variable). \\
$f_{p,m}^{(i,j)}$ & Traffic flow of mode $m$ on path $p$ between OD pair $(i,j)$, (Integer variable). \\
$f_{p,CD}^{(i,j)r,s}$ & Traffic flow of carpooling drivers on path $p$ between $(i,j)$ stopping at nodes $r$ and $s$, (Integer variable). \\ 
$y_{p,t}^{(i,j)}$ & Flow of empty or occupied service vehicles on path $p$ between OD pair $(i,j)$ ;  $t \in \{e, EH, RS\}$, (Integer variable). \\
$y_{p,RS}^{i,r,s,j}$ & Flow of ridesharing occupied vehicles on path $p$ stopping at nodes $i,r,s$ and $j$, (Integer variable).  \\
$c_{p,m}^{(i,j)}$ & Generalized cost of path $p$ with travel mode $m$ between OD pair $(i,j)$. 
\end{tabular} \label{tab:2}} 
\end{table}

\subsection{Multi-Modal Transportation Network without Intermodality}
\label{subsect:withoutIntermodality}
Let us consider an urban transportation network represented as a directed graph $G(E, A)$, wherein links ($A$) represent physical road sections, and nodes ($E$) can represent intersections or zones, depending on the level of network aggregation. 
Let $q^{(i,j)}$ be the travel demand between an Origin-Destination (OD) pair $(i,j)$, $\forall i,j \in E$. This study assumes that the demand for each OD pair is given. This assumption is valid when analyzing scenarios within a specific time period, such as peak hours, during which travel demand can be accurately estimated and remains relatively stable. Commuters traveling between $(i,j)$ must simultaneously choose their mode and path. For travel modes, as illustrated in \newref{Figure}{fig1}, commuters can take their car and be solo drivers, be a bus or metro (M) passenger, carpooling driver (CD), or carpooling passenger (CP) as long as the matching takes place, e-hailing passenger (EH) or ridesharing passenger (RS) and thus, be matched with a service vehicle. Additionally, they can reach their destination by walking (W) or biking (B). Accordingly, we define $\Psi$ as the set of travel modes, $\Psi = \{car, bus, M, W, B, CD, CP, EH, RS\} $. 

\begin{figure}[h!]
 \centering
 \includegraphics[width=0.8\textwidth]{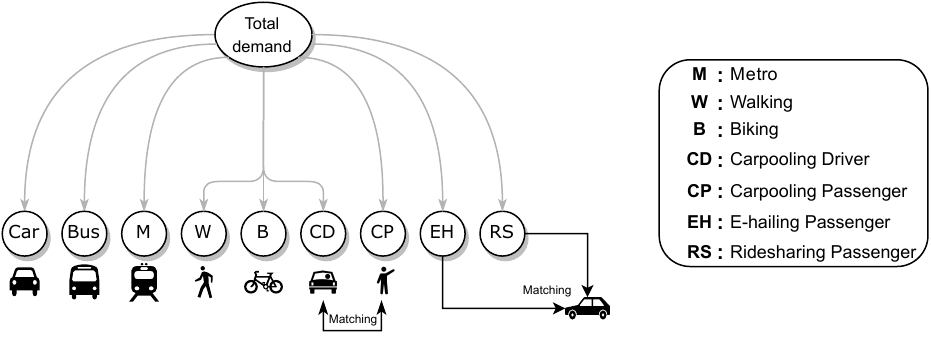}
 \caption{Available travel options for multi-modal networks, without intermodality.}
 \label{fig1}
\end{figure}
 
To consider interactions between the modes in terms of congestion effects, we represent the physical network by four subnetworks: road (RN), walking (WN), biking (BN), and metro (MN) network. Buses and cars both use the road network and contribute to traffic congestion. $N$ denotes the set of networks, $N = \{ RN, WN, BN, MN \}$. The model is formulated for a given time period to represent the state of the network, following two well-known principles: UE and SO. \cite{mirzahossein2021realistic} showed the validity of the TAP for a given time period compared with the reality represented by location-based data. Thus, the model is presented in a time-independent context, which keeps it tractable and analytically solvable.

\subsubsection{UE and SO Objective Function}
\label{sec:UESO}
This subsection presents the objective functions of UE and SO and how each component of these functions is calculated. In particular, we present the cost functions for commuters at the link level for all the modes mentioned in \newref{Figure}{fig1}.

Inspired by the formulation of travel costs used in the literature (see \newref{Table}{tab:1}), we define the generalized cost as the normalized sum of travel, waiting, service, and monetary costs. Additionally, through the Beckmann formulation~\citep{beckmann1956studies}, we solve the link-based user equilibrium while introducing path-based decision variables to consider mode and path choice for commuters. This results in a non-linear objective function as a classical form of resource allocation problem \citep{patriksson2008survey}. Moreover, when the flow variables are considered as integers, the problem results in a non-convex Mixed-Integer Bilinear Program (MIBLP), with the objective function defined as follows: 
\\
\begin{equation} \label{eq1}
min \; Z_{UE} = \sum_{n \in N} \sum_{a \in A_n}  \int_0^{x_a}  \alpha \cdot t_{a,n}(\omega) \; d\omega \quad + \quad  \sum_{a \in A} \sum_{m \in \Psi} [ \alpha (WT_{m,a} + ST_{m,a}) + C_{m,a} ] \cdot x_{a,m} 
\end{equation}
\\
\noindent This objective function represents the travel costs as the sum of each link's cost, which provides a general formulation that then handles intermodality in a straightforward way. In what follows, we present each component of this equation, defined separately for each mode. In equation (\ref{eq1}), $\alpha$ denotes the value of time, $x_{a,m}$ is the flow of mode $m$ on link $a$ and $x_a$ denotes the total flow of link $a$. $t_{a,n}$ denotes the travel time on link $a$ (i.e., the in-vehicle time), which is calculated by Equation (\ref{eq2}) based on the network to which the link belongs. For the road network (RN), the link travel time increases with the flow. In this study, we use the well-known Bureau of Public Roads (BPR) function~\citep{manual1964urban}. {Note that other travel time functions, such as the queue volume delay function (QVDF), proposed by \cite{zhou2022meso}, can also be considered.} For the other networks, the link travel time depends on the link length ($L_a$) and the average speed ($Sp_n$) of walking, biking, and metro, respectively. 

\begin{align} \label{eq2}
\nonumber\\
     t_{a,n} (x) = 
     \begin{cases}
           BPR (x) \, \text{or} \, QVDF (x)  & if  \; n = RN \\
           \frac{L_a}{Sp_n} & if \; n \in  \{MN, WN, BN\} \\ 
     \end{cases}  \\
\nonumber
\end{align}

Equation (\ref{eq3}) defines $WT_{m,a}$ as the waiting time of mode $m$ on link \textit{a}. A commuter has a waiting time at the trip's start, so it is only counted for links attached to an origin node (ensured by $o_a$,  {a binary coefficient which equals 1 if link $a$ is an origin link; 0 otherwise}). For PT, the waiting time is assumed to be half of the headway {(defined in a time-independent context as the inverse of frequency which represents the number of vehicles during the considered time period for the model)}. For SMSs, based on \cite{nourinejad_ride-sourcing_2020}, we assume a system-wide value for the waiting time, which depends on the number of passengers in the network ($q_m$) and an exogenous meeting rate ($R_m$). Since the number of passengers in the system is a decision variable, the waiting time for a shared mobility service is also variable and will proportionally increase with the demand. $R_m$ expresses the rate of pickups in the network and can be estimated by the Cobb-Douglas function~\citep{wei_modeling_2020}. Here, we set an average constant rate that should be estimated from real-world SMSs data. It should be noted that this value is only used to estimate the waiting time for SMSs and does not reflect the matching mechanism, which is handled by the decision variables of the model, as explained in the next section.

\begin{align} \label{eq3}
\nonumber\\
     WT_{m,a} = 
     \begin{cases}
           0  & if  \; m \in \{car, CD, W, B\} \\
           \frac{1}{2freq_{m,a}} \cdot o_a & if \;m \in \{Bus, M\} \\
           \frac{q_m}{R_m} \cdot o_a & if \;m \in \{CP,RS,EH\} \\
     \end{cases} 
\end{align}
\\
$ST_{m,a}$ in Equation (\ref{eq4}) is the service time of mode $m$ on link $a$. For PT, it denotes the average time spent at each stop ($S_{Bus}$, $S_{M}$ for bus and metro). For car or bike users, it represents the parking time ($P_{car}$ and $P_B$) and is only counted for destination links (represented by $d_a$,  {a binary coefficient which equals 1 if link $a$ is a destination link; 0 otherwise}). $ST_{m,a}$ for SMSs represent the time to be matched with a driver/passenger and the boarding and drop-off times ($S_m$; $m \in \{CP, RS, EH\}$) on origin and destination links. It also implicitly includes any additional cost due to SMSs use and detours. Additionally, a CD experiences parking time at the destination. Walking does not have a service time.
\\
\begin{equation} \label{eq4}
     ST_{m,a} = 
     \begin{cases}
           0 & if \; m =W \\
           S_m & if \; m \in \{Bus, M\} \\
           S_m \cdot (o_a + d_a) & if \; m \in \{CP, RS, EH\} \\
           S_m \cdot (o_a + d_a) + P_m \cdot d_a & if \; m = CD \\
           P_m \cdot d_a  & if \; m \in \{ car, B\} 
     \end{cases} 
\end{equation}
\\
Equation (\ref{eq5}) denotes the monetary cost of mode $m$ on link \textit{a}. Based on~\cite{du_modeling_2022}, for solo drivers, it includes the cost of using the car (i.e., fuel cost) for the traveled distance and the parking fare ($PF_{car}$) at the destination. For all other passengers (PT and SMSs passengers), $C_{m,a}$ represents the trip fare to be paid ($TF_{m,a}$). For carpooling drivers, the trip fare received from passengers is retrieved from the cost of using their car. For walking and biking, we assume no monetary cost.

\begin{align} \label{eq5}
\nonumber\\
     C_{m,a} = 
     \begin{cases}
           \gamma \cdot L_a + PF_m \cdot d_a & if  \; m = car \\
           TF_{m,a}  & if \; m \in \{ Bus, M, CP, EH, RS\} \\
           \gamma\cdot L_a + PF_m \cdot d_a - TF_{m,a} & if  \; m = CD \\
           0 & if \; m \in \{W,B\} 
     \end{cases} \\
\nonumber     
\end{align}

\noindent Using these same definitions, we define the SO objective function to minimize the system's total cost as:

\begin{align} \label{eq6}
\nonumber\\
min \; Z_{SO} = \sum_{n \in N} \sum_{a \in A_n}  \alpha \cdot t_{a,n}(x_{a})\cdot x_{a} \quad + \quad  \sum_{a \in A} \sum_{m \in \Psi} [ \alpha (WT_{m,a} + ST_{m,a}) + C_{m,a} ] \cdot x_{a,m}  \\
\nonumber 
\end{align}

It is worth mentioning that the considered mathematical program, for both UE and SO, is bilinear because of the second term of the objective function, in which the waiting time for shared mobility services is multiplied by the modal link flow $x_{a,m}$. As for the BPR function in the first term of $Z_{UE}$ and $Z_{SO}$, we apply a linearization procedure as described in \ref{app::linear}

In what follows, we present the model's constraints before introducing intermodality.

\subsubsection{Model's Linear Constraints for Mode and Path Choice}

Traffic flows, particularly for SMSs, often involve discrete units (individual or vehicle units), which makes an integer representation more suitable for accurately capturing the traffic dynamics, particularly the interactions between modes and the matching process between drivers and passengers for SMSs.
Following this, we define the flow variables of the model as integers, aligning with the trip-based approach commonly used in the literature (e.g, \cite{geurs_multi-modal_2016}, \cite{liu_doubly_2017}, \cite{narayan_does_2019}). However, it is worth mentioning that our model is still valid with continuous flow variables. In what follows, we will progressively introduce these variables and the model's constraints, ensuring the solution's coherence with reality.  

The objective functions previously described are formulated in terms of link flow and cost. However, commuters are confronted with mode and path choices. Thus, let $f_{p,m}^{(i,j)}$ be the path flow variable. Equation (\ref{eq25}) presents the link-path flow conversion, and Equation (\ref{eq7}) represents the demand conservation constraint.

\begin{align}
\nonumber\\
    &x_{a,m} = \sum_{i,j \in E} \sum_{p \in P_{ij}^m} f_{p,m}^{(i,j)} \cdot \delta_{a,p,m}^{(i,j)} & \forall a \in A \quad ; \; \forall m \in \Psi \label{eq25} \\
    &q^{(i,j)} = \sum_{m \in \Psi} \sum_{p \in P_{ij}^m} f_{p,m}^{(i,j)}  & \forall i,j \in E \label{eq7} 
\end{align}
\\
{Note that $\delta_{a,p,m}^{(i,j)}$ in Equation \ref{eq25} can be aligned with the definition in \cite{larsson1995augmented}, which can be used to define a capacitated traffic assignment model. However, the objective function may lose strict convexity due to the multimodal definition at the path level. Consequently, the Lagrangian multiplier cannot be explicitly derived for the path flow capacity constraints \citep{correa2004selfish}. On the other hand, similar to \cite{mahmoudi2021many}, we can track the mode choice of all users along their paths, allowing us to develop a detailed path-based model that includes intermodal options in addition to link-based traffic assignment.} Here, we present the constraints for every mobility option considered. Equations (\ref{eq8})-(\ref{eq10}) are for the carpooling service. We define $f_{p,CD}^{(i,j)\;r,s}$ as the number of carpooling drivers traveling between $(i,j)$, with path $p$ and stopping at nodes $r$ and $s$ to pick up or drop off their passengers. Note that a CD can have more than one passenger onboard. Thus, this variable represents the matching process for carpooling in terms of OD pairs. Consequently, $f_{p,CD}^{(i,j)\;i,j}$ expresses the number of CD having passengers with the same OD pair $(i,j)$. Constraint (\ref{eq8}) ensures that the total number of CD between $(i,j)$ is the sum of the drivers stopping at the different passengers' OD pairs on path $p$. 
To formulate the correspondence between the drivers' and passengers' paths, let us define the following sets: 
\begin{itemize}
    \item $Q_p^{(r,s)} $: Set of paths $l$ between $(r,s)$, where the path $p$ is a sub-path of $l$. 
    \item $R_p^{(r,s)}$: Set of paths $l$ $(r,s)$, where $l$ is a sub-path of $p$.
\end{itemize} \vspace{0.3cm}

Constraint (\ref{eq9}) ensures drivers can pick up all carpooling passengers. Recall that the drivers may or may not have the same OD pair as the passengers. \newref{Figure}{figPSP} illustrates the correspondence between the carpooling passenger and driver paths, expressed by these constraints. Constraint (\ref{eq10}) ensures that all CD stopping at nodes $r$ and $s$ have enough passengers between $(r,s)$. We formulated these constraints regarding path flow variables to ensure the passenger-driver conservation over the path $p$. 

\begin{figure}[h!]
 \centering
 \includegraphics[width=0.5\textwidth]{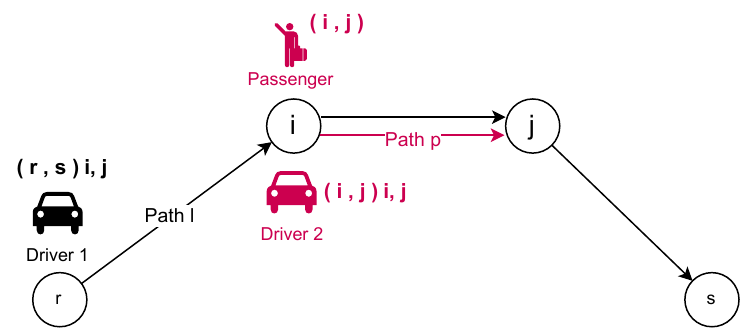}
 \caption{Correspondence between the carpooling passenger's and driver's path.}
 \label{figPSP}
\end{figure}
\begin{align}
    &f_{p,CD}^{(i,j)} = \sum_{r \in E} f_{p,CD}^{(i,j)\; r} & \forall i,j \in E \quad ; \; \forall p \in P_{ij}^{CD} \label{eq8} \\
    & f_{p,CP}^{(i,j)} \leq CAP_{CP} \cdot [\sum_{r,s \in E} \; \sum_{l \in Q_p^{(r,s)} } f_{l,CD}^{(r,s)\; i,j} ] &  \forall i,j \in E \quad ;\; \forall p \in P_{ij}^{CP}\label{eq9} \\
    &f_{l,CD}^{(i,j)\; r,s} \leq f_{p,CP}^{(r,s)} &  \forall r,s,i,j \in E \quad ;\; \forall l \in P_{ij}^{CD} \quad ;\; \forall p \in R_l^{(r,s)} \label{eq10} \\
\nonumber 
\end{align}

For e-hailing, $y_{p,EH}^{(i,j)}$ denotes the number of fleet vehicles having an e-hailing passenger between $(i,j)$ with path $p$. Constraint (\ref{eq11}) ensures that this variable equals the corresponding e-hailing path flow variable. 

\begin{align}
    &f_{p,EH}^{(i,j)} = y_{p,EH}^{(i,j)} & \forall i,j \in E \quad ; \; \forall p \in P_{ij}^{EH} \label{eq11} \\
\nonumber
\end{align}

Here, we propose a ride-sharing service formulation. Let $y_{p,RS}^{i,r,s, j}$ represent the number of fleet vehicles participating in a ridesharing service and stopping at nodes $i$,$j$,$r$ and $s$ to pick up or drop off passengers. The ordering of the nodes in this variable depends on whether a First In First Out (FIFO) or Last In First Out (LIFO) service is considered, as illustrated in \newref{Figure}{fig2}. 

\begin{figure}[h!]
 \centering
 \includegraphics[width=0.8\textwidth]{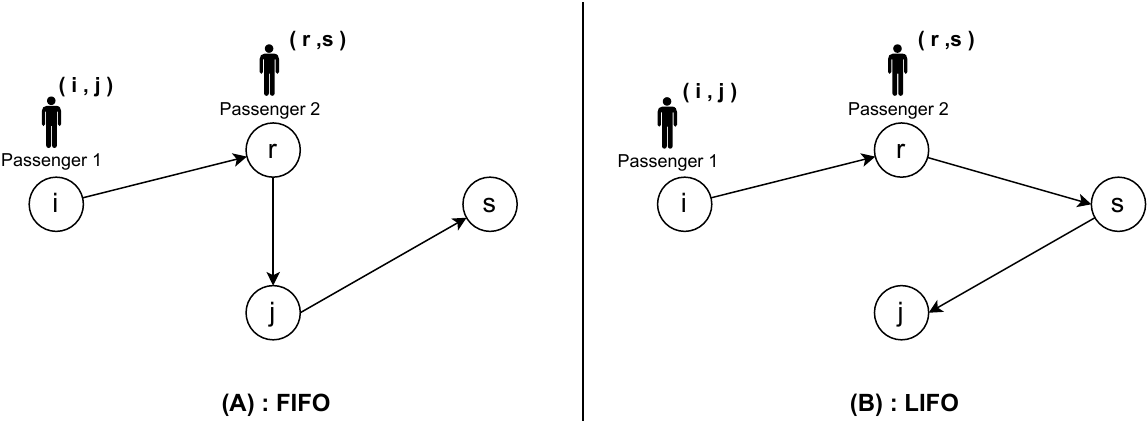}
 \caption{Order of ridesharing service. (A): The first passenger picked up is the first dropped off; (B): The first passenger picked up is the last dropped off.  }
 \label{fig2}
\end{figure}

Equation (\ref{eq12}) defines the number of ridesharing vehicles traveling between $(i,j)$ with path $p$. Similar to the carpooling service, Constraint (\ref{eq13}) ensures enough seats for the passengers, and Constraint (\ref{eq14}) ensures enough passengers for the ridesharing cars. $CAP_{RS}$ denotes the ridesharing passenger capacity, which is the same for all RS vehicles.

\begin{align}
    &y_{p,RS}^{(i,j)} = \sum_{r,s\in E} y_{p,RS}^{i,r,s,j}&  \forall i,j \in E \quad ; \; \forall p \in P_{ij}^{RS} \label{eq12}\\ 
    & f_{p,RS}^{(i,j)} \leq CAP_{RS} \; \cdot [\sum_{r,s \in E} \; \sum_{l \in Q_p^{(r,s)} } (y_{l,RS}^{r,i,j, s} + y_{l,RS}^{r,i,s, j} ) \; ]&  \forall i,j \in E \quad ; \; \forall p \in P_{ij}^{RS}\label{eq13} \\
    & y_{l,RS}^{i,r,s, j} +  y_{l,RS}^{i,r,j, s}\leq f_{p,RS}^{(r,s)} &  \forall r,s,i,j \in E \quad ;\; \forall l \in P_{ij}^{RS} \quad ; \; \forall p \in R_p^{(r,s)} \label{eq14} \\
    \nonumber
\end{align}

Constraint (\ref{eq15}) ensures that the ridesharing passengers are sharing the trip with another passenger (not an e-hailing service). Note that $y_{p,RS}^{i,i,j, j}$ represents the ridesharing vehicles on path $p$ having passengers with the same OD pair $(i,j)$. This is different from carpooling, as the driver providing the service is not a commuter himself.

\begin{align}
\nonumber \\
    &f_{p,RS}^{(i,j)} \geq 2 \cdot y_{p,RS}^{i,i,j,j} + \sum_{r,s \in E} \sum_{l \in Q_p^{(r,s)}} (y_{l,RS}^{r,i,j,s} + y_{l,RS}^{r,i,s,j}) & \forall i,j \in E \quad ; \; \forall p \in P_{ij}^{RS} \label{eq15}
\end{align} 
\\
Constraints (\ref{eq16}) and (\ref{eq17}) define the number of occupied fleet cars as the sum of cars participating in the e-hailing and ridesharing services. 

\begin{align}
\nonumber \\
    &q_o^{(i,j)\;r} = \sum_{p \in P_{ij}^{RS}}  y_{p,RS}^{i,r,s,j} + y_{p,RS}^{i,s,r,j} & \forall i,j,r \in E \label{eq16} \\
    &q_o^{(i,j)} =  \sum_{r \in E} q_{o}^{(i,j)\;r} + \sum_{p \in P_{ij}^{EH}} y_{p,EH}^{(i,j)} & \forall i,j \in E \label{eq17} 
\end{align}

Based on the formulation of the SMSs, we need to conserve the passenger and driver flow over the paths. $y_{p,e}^{(i,j)}$ denotes the number of empty fleet vehicles traveling between $(i,j)$ with path $p$.  
Constraint (\ref{eq19}) ensures, for every destination node, that the number of entering vehicles equals the number of exiting vehicles, {while keeping the number of occupied and empty vehicles coherent, based on the pickups (i.e. origins) and drop-offs (i.e. destinations)}. Constraint (\ref{eq20}) ensures, for every origin node, that there are enough empty cars to pick up the e-hailing and ridesharing passengers. Constraint (\ref{eq21}) ensures the fleet size ($|V|$) conservation. Constraints (\ref{eq22})-(\ref{eq24}) are integrality conditions on the decision variables.

\begin{align}
    &q_e^{(i,j)} = \sum_{p \in P_{ij}^{EH}} y_{p,e}^{(i,j)} & \forall i,j \in E \label{eq18}\\
    &\sum_{i,j \in E} q_o^{(i,j)\;s} + \sum_{j \in E} q_e^{(j,s)} = \sum_{j \in E} q_o^{(s,j)} + \sum_{j \in E} q_e^{(s,j)} & \forall s \in D \label{eq19} \\
    &\sum_{i \in E} q_e^{(i,j)} \geq \sum_{s \in E} q_o^{(j,s)} & \forall j \in O \label{eq20} \\
    &\sum_{i,j \in E} q_o^{(i,j)} + q_e^{(i,j)} = |V| \label{eq21} \\
    &f_{p,m}^{(i,j)}, f_{p,CD}^{(i,j)r,s} \geq 0 & \forall p \in P_{ij}^m \quad ; \; \forall i,j,r,s \in E \quad ; \; \forall m \in \Psi \label{eq22} \\ 
    &y_{p,m}^{(i,j)} \geq 0 & \forall p \in P_{ij}^m \quad ; \; \forall i,j \in E \quad ; \; \forall m \in \{EH, RS, e\} \label{eq23} \\ 
    &y_{p,RS}^{i,r,s,j} \geq 0 & \forall p \in P_{ij}^{RS} \quad ; \; \forall i,j,r,s \in E \label{eq24} 
\end{align}

Although we do not explicitly consider maximum detour constraints in the model, we ensure for SMSs that the detour remains acceptable by the passengers. This is done through the pricing scheme. In particular, the fare to be paid by the passengers increases with the length of the trip ($TF_{m,a}$ is specified for each link, and the more links are traversed, the higher will be the price, and potentially travel time). Thus, the best compromise between the detour due to sharing, and the monetary gains from sharing will be found by the system at the optimal solution.

\subsubsection{Link Flow and Path Flow Correspondence} \label{sec:linkflow}
The UE and SO objective functions, formulated in \newref{Section}{sec:UESO}, require the calculation of $x_a$ and $x_{a,m}$. Following the Beckmann formulation, $x_{a,m}$ is defined by Equation (\ref{eq25}). However, to calculate the $x_a$ variable used in links travel time, we consider the travel modes contributing to traffic congestion, occupied and empty fleet vehicles, and PT flow. Considering PT flow means we need to count PT units independently of the passenger flow. This is consistent with reality since PT adheres to a well-defined schedule and operates even with no passengers. It is worth mentioning that many studies (e.g.,~\cite{di_unified_2019},~\cite{du_modeling_2022}, and~\cite{zhou_multimodal_2022}) formulate only non-empty PT flow, which underestimates the impact of PT on congestion levels. In our study, and since we consider a time-independent context, we use the frequency of PT on every link to count the units, as described by Equation (\ref{eq26}). 

\begin{align} \label{eq26}
\nonumber \\
    x_a = \sum_{m \in \Psi'} x_{a,m} + \sum_{l \in PT} freq_{l,a} + \sum_{i,j \in E} \sum_{p \in P_{ij}^{RS}} y_{p,RS}^{(i,j)} \cdot & \delta_{a,p,RS}^{(i,j)}  + \sum_{i,j \in E} \sum_{p \in P_{ij}^{EH}} (y_{p,EH}^{(i,j)} + y_{p,e}^{(i,j)} ) \cdot \delta_{a,p,EH}^{(i,j)} \\  & \forall a \in A \quad ; \; \Psi' = \{car, CD, W, B\} \quad ; \; PT = \{bus, M\} \nonumber 
\end{align}
\\
Regarding the carpooling service, we only count the driver flow. For the ridesharing service, we count the number of cars instead of the passenger flow since a vehicle may carry several passengers. This enables us to keep the link and path flows as the number of vehicles, aligning with an integer value. All modes are defined at this point, and their travel cost functions are formulated. Now, we can introduce intermodality.

\subsection{Multi-Modal Transportation Network with Intermodality}
\label{subsect:withIntermodality}

This section introduces the transformation of the model's constraints and variables to consider intermodality. In this study, an inter-mode is any combination of the previously enumerated modes (\newref{Figure}{fig1}) and is considered a separate mode. For example, Park-and-Ride as a travel option is the commuter who uses the car to a transfer node and then takes the metro~\citep{murphy_shared_2016}.

The proposed model can consider any logical combination of modes as a travel option without further modifications. Thus, the set of modes $\Psi$ is extended to include all modes in the form of $I_{m1,m2}$ where $m1$ and $m2 \in \Psi$. Let us consider a trip with \textit{\textbf{k}} modes. This intermodal trip can be represented by a sequence of modes (simple or intermodes) within the set $\Psi$, as illustrated by \newref{Figure}{fig3_1}. For instance, if the commuter uses the car to reach a metro station (park and ride mode), then a ride-sharing service for the last mile, this trip can be modeled as $I_{m1,RS}$ with $m1 = I_{car,M}$. However, as a first step in analyzing travel behavior with intermodality, the experiments conducted in this study consider the inter-modes represented in \newref{Figure}{fig3}. 

\begin{figure}[h!]
 \centering
 \includegraphics[width=0.3\textwidth]{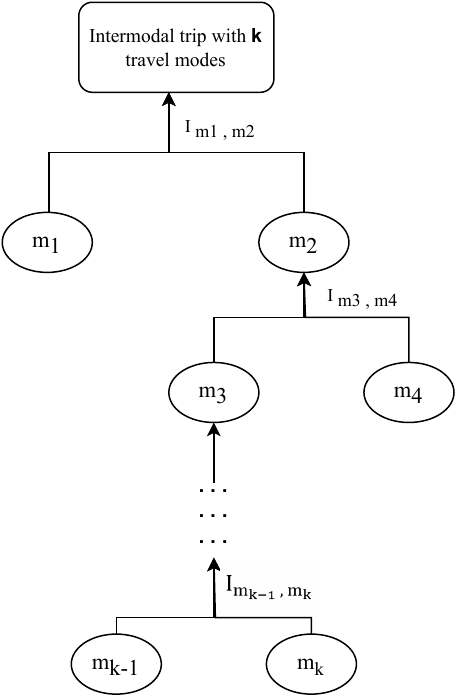}
 \caption{Representation of an intermodal trip with \textbf{\textit{k}} modes. }
 \label{fig3_1}
\end{figure}

As we are considering every inter-mode as a separate mode, commuters between $(i,j)$ with an inter-mode $m = I_{m1,m2}$ on path $p$ ($p = p1 \oplus p2$, where $\oplus$ is the concatenation of two paths), use mode $m1$ with path $p1$, up to a transfer node. Then, they use $m2$ with path $p2$ to their destination. The number of commuters between $(i,j)$ with an inter-mode $m$ using path $p$ is represented by the path flow variable $f_{p,m}^{(i,j)}$. Note that the model can consider a combination of more than two modes $I_{m1,m2}$, where $m1$ is a 'simple' mode and $m2$ is an inter-mode (and vice-versa).

\begin{figure}[h!]
 \centering
 \includegraphics[width=0.8\textwidth]{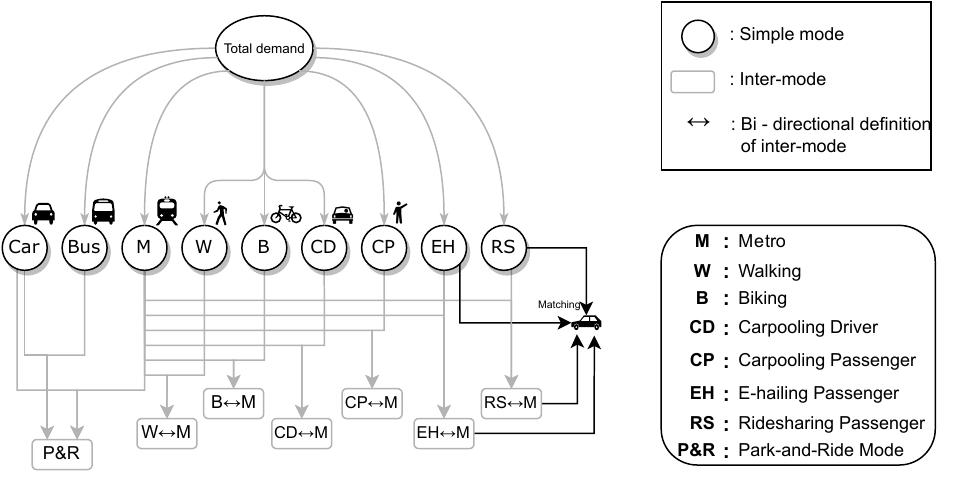}
 \caption{Available travel options for multi-modal networks, with intermodality. }
 \label{fig3}
\end{figure}

Let $PI( (i,j), s)$ be the set of all inter-modal paths between $(i,j)$ in which the node $s$ is the transfer node. Constraint (\ref{eq27}) ensures that the number of commuters using their car as part of an inter-mode and transferring at node $s$ does not exceed the parking capacity at node $s$, represented by $Pk\_cap_s$. 

\begin{align}
\nonumber \\
    &\sum_{i,j \in E} \sum_{p \in PI( (i,j),s)}  f_{p,I_{car,M}}^{(i,j)} + f_{p,I_{CD,M}}^{(i,j)} \leq Pk\_cap_s & \forall s \in TR  \label{eq27} 
\end{align}
\\
To clearly represent intermodality, let I(m) denote an inter-mode including the simple mode $m$ (either as the first part or the second part of the inter-mode). m(p) denotes the part of path $p$ where the mode $m$ is used, and TR(p) denotes the transfer node of the inter-modal path $p$. To consider intermodality with carpooling, we add Constraints (\ref{eq28}) and (\ref{eq29}) to the model. Constraint (\ref{eq28}) ensures that drivers can pick up the passengers participating in carpooling as part of an inter-mode. The right-hand side of the constraint considers drivers participating in carpooling as part of an intermodal trip, as well as carpooling drivers between any OD pair $(u,v)$ passing by the transfer node $t$ of the passengers. This is illustrated by \newref{Figure}{figIPSP}. In other words, a passenger can be picked up by a driver going from the same origin to the same destination or by a driver going from the origin to the transfer node $t$.

\begin{figure}[h!]
 \centering
 \includegraphics[width=0.5\textwidth]{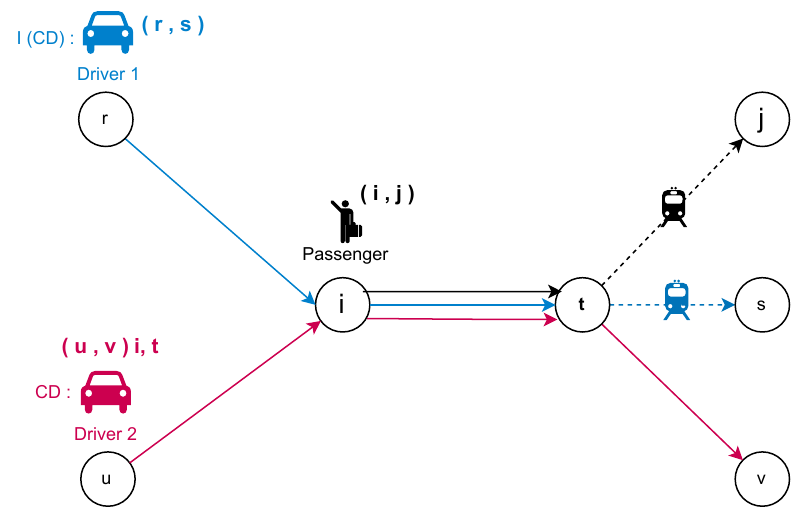}
 \caption{Correspondence between a carpooling passenger and drivers' paths with intermodality.}
 \label{figIPSP}
\end{figure} 

{\footnotesize
\begin{align} \label{eq28} 
    f_{p,I(CP)}^{(i,j)} \leq CAP_{CP} \cdot (\sum_{I(CD) \in \Psi } \sum_{r,s \in E} \sum_{k \in \{ x \in P_{rs}^{I(CD)} \; ; \; CP(x)=CP(p)\} } &f_{k,I(CD)}^{(r,s)} + \sum_{u,v \in E} \sum_{l \in P_{uv}^{CD}} f_{l,CD}^{(u,v)i,t} + f_{l,CD}^{(u,v)t,j} )  \\
    &  \forall i,j \in E \quad ; \; \forall p \in P_{ij}^{I(CP)} \quad ; \; t = TR(p)  \nonumber
\end{align}}

Constraint (\ref{eq29}) ensures that all carpooling drivers between $(i,j)$ have enough carpooling passengers to pick up. This means that no CD is driving without a passenger on board, either for the full trip (i.e., same OD pair as the passenger) or part of the trip.

{\footnotesize
\begin{align} \label{eq29} 
\nonumber \\
    f_{p,I(CD)}^{(i,j)} \leq \sum_{I(CP) \in \Psi } \sum_{r,s \in E} \sum_{p' \in \{ x \in P_{rs}^{I(CP)} \; ; \; CP(x)=CP(p)\} }  f_{p',I(CP)}^{(r,s)} + & \sum_{u,v \in E} \sum_{k \in \{ x \in P_{uv}^{CP} \; ; \; x=CP(p)\} }  f_{k,CP}^{(u,v)} \quad  \forall i,j \in E \quad ; \; \forall p \in P_{ij}^{I(CD)} 
\end{align}}

\noindent Thus, we need to reformulate Constraints (\ref{eq9}) and (\ref{eq10}) to include intermodality.

{\footnotesize
\begin{align}\label{eq30}
\nonumber \\
    f_{p,CP}^{(i,j)} \leq CAP_{CP} \cdot (\sum_{r,s \in E} \; \sum_{k \in Q_p^{(r,s)} } f_{k,CD}^{(r,s)\; i,j} + \sum_{I(CD) \in \Psi } \sum_{u,v \in E} & \sum_{k \in \{ x \in P_{uv}^{I(CD)} \; ; \; CP(x)=p\} }  f_{k,I(CD)}^{(u,v)} )   \quad \; \forall i,j \in E \quad ;\; \forall p \in P_{ij}^{CP}  
\end{align}} 

{\footnotesize
\begin{align}\label{eq31}
    f_{p,CD}^{(i,j)\; r,s} \leq f_{k,CP}^{(r,s)} + \sum_{I(CP) \in \Psi } \sum_{u,v \in E} \sum_{k \in \{ x \in P_{uv}^{I(CP)} \; ; \; CP(x)=p\} }  f_{k,I(CP)}^{(u,v)} &  \quad  \forall r,s,i,j \in E \quad ;\; \forall p \in P_{ij}^{CD} \quad ;\; \forall k \in R_p^{(r,s)}   
\end{align}}

\noindent Similarly, for e-hailing and ridesharing, Constraint (\ref{eq32}) is added to the model, and Constraints (\ref{eq11}) and (\ref{eq14}) are updated accordingly.

{\footnotesize
\begin{align} \label{eq32} 
\nonumber \\
    f_{p,I(RS)}^{(i,j)} \leq CAP_{RS} \cdot  \sum_{r,s \in E} \sum_{k \in Q_p^{(r,s)}} (y_{k,RS}^{r,i,TR(p),s} + y_{k,RS}^{r,i,s,TR(p)} + & y_{k,RS}^{r,TR(p),j,s} +  y_{k,RS}^{r,TR(p),s,j} )  \quad  \forall i,j \in E \quad ; \; \forall p \in P_{ij}^{I(RS)} 
\end{align}}

\begin{align}
    &f_{p,EH}^{(i,j)} + \sum_{I(EH) \in \Psi } \sum_{u,v \in E} \sum_{k \in \{ x \in P_{uv}^{I(EH)} \; ; \; EH(x)=p\} }  f_{k,I(EH)}^{(u,v)}  = y_{p,EH}^{(i,j)} & \forall i,j \in E \quad ; \; \forall p \in P_{ij}^{EH} \label{eq33}
\end{align}

{\footnotesize
\begin{align} \label{eq34} 
    y_{p,RS}^{i,r,s, j} +  y_{p,RS}^{i,r,j, s} \leq f_{k,RS}^{(r,s)} + \sum_{I(RS) \in \Psi } \sum_{u,v \in E} \sum_{k \in \{ x \in P_{uv}^{I(RS)} \; ; \; RS(x)=p\} } & f_{k,I(RS)}^{(u,v)} \quad \forall r,s,i,j \in E \quad ;\; \forall p \in P_{ij}^{RS} \quad ; \; \forall k \in R_p^{(r,s)}  
\end{align}}

\noindent To have the link-path flow correspondence, we update the $x_{a,m}$ variable as described by Equation (\ref{eq35}).

\begin{align} \label{eq35}
\nonumber \\
    x_{a,m} = \sum_{i,j \in E} \sum_{p \in P_{ij}^m}  f_{p,m}^{(i,j)}\cdot \delta_{a,p,m}^{(i,j)} + \sum_{I(m) \in \Psi } \sum_{u,v \in E} \sum_{k \in P_{uv}^{I(m)}} f_{k,I(m)}^{(u,v)}  \cdot \delta_{a,m(k),m}^{(u,TR(k))} \qquad \forall a \in A \quad ; \; \forall m \in \Psi 
\end{align}
\\
\noindent Thus, the final comprehensive multi-modal mathematical formulation of UE (UE-CMF) is: 
{\footnotesize
\begin{gather*}
        \tag{UE-CMF}
        \label{P1}
        min \; Z_{UE} = \sum_{n \in N} \sum_{a \in A_n}  \int_0^{x_a}  \alpha \cdot t_{a,n}(\omega) \; d\omega \quad + \quad  \sum_{a \in A} \sum_{m \in \Psi} [ \alpha (WT_{m,a} + ST_{m,a}) + C_{m,a} ] \cdot x_{a,m} \\ 
        \notag
        \quad \textit{s.t.} \\  \qquad \left\{\begin{array}{l @{\quad} l r}
        \text{Cost functions} & & (\ref{eq2})-(\ref{eq5}) \\
        x_{a,m} = \sum_{i,j \in E} \sum_{p \in P_{ij}^m} f_{p,m}^{(i,j)} \cdot \delta_{a,p,m}^{(i,j)} & \forall a \in A  ; \; \forall m \in \Psi & (7) \\
         q^{(i,j)} = \sum_{m \in \Psi} \sum_{p \in P_{ij}^m} f_{p,m}^{(i,j)}  & \forall i,j \in E & (8)\\
        f_{p,CD}^{(i,j)} = \sum_{r \in E} f_{p,CD}^{(i,j)\; r} & \forall i,j \in E  ; \; \forall p \in P_{ij}^{CD} & (9)\\
        y_{p,RS}^{(i,j)} = \sum_{r\in E} y_{p,RS}^{(i,j)\;r}&  \forall i,j \in E  ; \; \forall p \in P_{ij}^{RS} & (13)\\ 
        f_{p,RS}^{(i,j)} \leq CAP_{RS} \; \cdot [\sum_{r,s \in E} \; \sum_{l \in Q_p^{(r,s)} } (y_{l,RS}^{r,i,j, s} + y_{l,RS}^{r,i,s, j} ) \; ]&  \forall i,j \in E  ; \; \forall p \in P_{ij}^{RS} & (14) \\
        y_{l,RS}^{i,r,s, j} +  y_{l,RS}^{i,r,j, s}\leq f_{p,RS}^{(r,s)} &  \forall r,s,i,j \in E  ;\; \forall l \in P_{ij}^{RS}  ; \; \forall p \in R_p^{(r,s)} & (15) \\
        f_{p,RS}^{(i,j)} \geq 2 \cdot y_{p,RS}^{i,i,j,j} + \sum_{r,s \in E} \sum_{l \in Q_p^{(r,s)}} (y_{l,RS}^{r,i,j,s} + y_{l,RS}^{r,i,s,j}) & \forall i,j \in E  ; \; \forall p \in P_{ij}^{RS} & (16) \\
        q_o^{(i,j)\;r} = \sum_{p \in P_{ij}^{RS}}  y_{p,RS}^{i,r,s,j} + y_{p,RS}^{i,s,r,j} & \forall i,j,r \in E & (17) \\
        q_o^{(i,j)} =  \sum_{r \in E} q_{o}^{(i,j)\;r} + \sum_{p \in P_{ij}^{EH}} y_{p,EH}^{(i,j)} & \forall i,j \in E & (18) \\
        q_e^{(i,j)} = \sum_{p \in P_{ij}^{e}} y_{p,e}^{(i,j)} & \forall i,j \in E & (19)\\
        \sum_{i,j \in E} q_o^{(i,j)\;s} + \sum_{j \in E} q_e^{(j,s)} = \sum_{j \in E} q_o^{(s,j)} + \sum_{j \in E} q_e^{(s,j)} & \forall s \in D & (20) \\
        \sum_{i \in E} q_e^{(i,j)} \geq \sum_{s \in E} q_o^{(j,s)} & \forall j \in O & (21) \\
        \sum_{i,j \in E} q_o^{(i,j)} + q_e^{(i,j)} = |V| & & (22) \\
        f_{p,m}^{(i,j)}, f_{p,m}^{(i,j)r,s} \geq 0 & \forall p \in P_{ij}^m  ; \; \forall i,j,r,s \in E  ; \; \forall m \in \Psi & (23) \\ 
        y_{p,m}^{(i,j)} \geq 0 & \forall p \in P_{ij}^m  ; \; \forall i,j \in E  ; \; \forall m \in \{EH, RS, e\} & (24) \\ 
        y_{p,RS}^{i,r,s,j} \geq 0 & \forall p \in P_{ij}^{RS} ; \; \forall i,j,r,s \in E & (25) \\
        \text{Intermodality Constraints}    &   & (\ref{eq26})-(\ref{eq35})  \\ 
        \end{array}\right.
\end{gather*}
}
To address the SO principle, only the objective function in \ref{P1} needs to be replaced by Equation (\ref{eq6}), and the problem can be treated as any minimization problem and solved by common (open-source or commercial) solvers. {Moreover, if we keep the same configuration of the model while adding an adequate stochastic term in the objective function, the program can be converted to an SUE calculation.}

In what follows, we discuss the theoretical aspects of the MIBLP \ref{P1}.

\paragraph{\textbf{Proposition 1}} Solution of optimization problem \ref{P1} satisfies user equilibrium conditions {with respect to a delay induced by the use of SMSs, and subject to capacity constraints}.
\subparagraph{\textbf{Proof.}} See \ref{app::UE}.

{It is worth mentioning that this redefinition of the UE principle is due to the system-wide formulation of the waiting time for SMSs. If we alternatively consider a link-based formulation for the waiting time, the program will be converted to the classical UE program, with capacity constraints.}

Since the objective function in \ref{P1} is non-convex, the uniqueness of equilibrium solutions can be proven at the link level if all cost functions remain strictly monotone~\citep{ameli_evolution_2022}. However, the uniqueness with respect to path flows is not straightforward. In other words, there could exist an infinite number of path flow solutions leading to the same unique link flow solution~\citep{rakha_traffic_2009}. This aspect is further emphasized with the existence of multiple travel modes that can generate the same link flows at equilibrium.

Additionally, the introduction of integer variables does not affect the UE conditions of the program, but may lead to a different solution compared to the continuous one. In \cite{rosenthal_network_1973}, the authors showed that considering the flow variables of a UE program as integers is equivalent to treating the problem as an n-person non-cooperative game to which the solution is a Nash equilibrium. However, considering the commuters as discrete entities rather than a continuous flow can lead to differences in their costs. To illustrate this, let us consider the scenario in \newref{Figure}{figInt_SC}, with a single unit of travel between the OD pair $(A,B)$. Two paths connect point $A$ to point $B$, and the cost of each path equals to its flow. For this example, we do not consider multiple travel modes since the mode choice can be seen as an extended path choice. Thus, we only focus on the choice of paths with the minimum cost to illustrate the equilibrium state.

\begin{figure}[h!]
 \centering
 \includegraphics[width=0.3\textwidth]{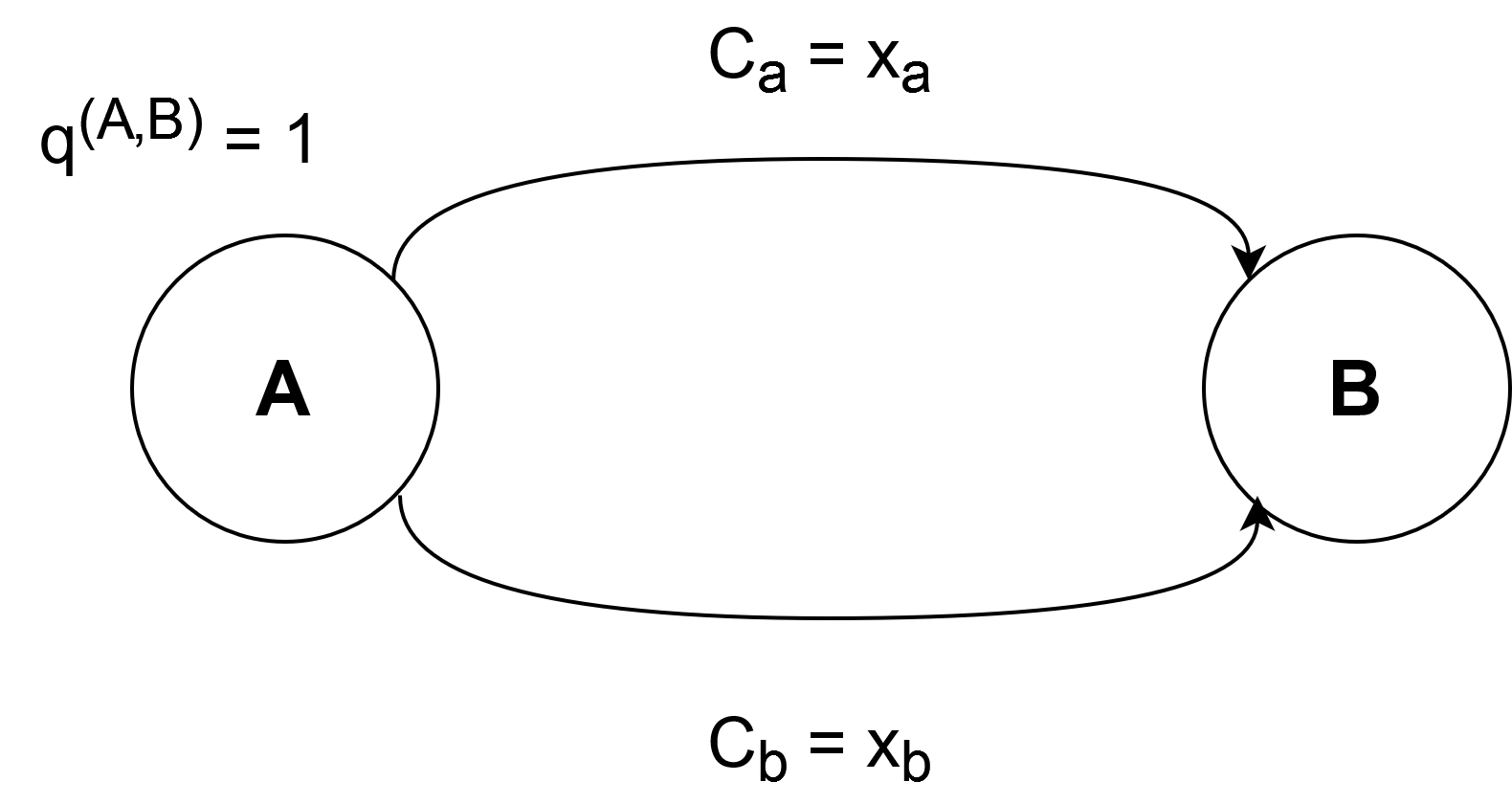}
 \caption{Example to illustrate the effect of integer flows.}
 \label{figInt_SC}
\end{figure} 

In this scenario, considering the travel unit as a continuous flow will result in a flow of $0.5$ on link a and $0.5$ on link b. However, we will have a flow of $1$ on link $a$ (or link $b$) if the travel unit is atomic, i.e., integer flow. An important point to raise is that the total cost will equal $0.5$ with a continuous setting, while it equals $1$ with an integer setting. This means that when we consider integer values, commuters can experience higher costs (delays). However, considering the continuous setting will not give us the opportunity to analyze commuters' choices since it does not provide an intuitive meaning for the flow. Consequently, in our experimental analysis, we consider the integer configuration in order to provide a meaningful and more realistic analysis.

\subsection{Price of Anarchy Formulation}
\label{subsect:PoA}

The Price of Anarchy (PoA) is a term in economics and game theory used to measure how the efficiency of a system degrades due to the selfish and conflicting behavior of its agents. In the context of transportation, this concept was first introduced by~\cite{koutsoupias1999worst}. The authors formulated the network routing problem as a non-cooperative congestion game. They proposed a measure of the performance degradation caused by a lack of coordination as the worst-case ratio of social welfare. This ratio is achieved by user equilibrium and by a socially optimal set of strategies~\citep{roughgarden_how_2002},~\citep{christodoulou_price_2005}. 

In this study, the Price of Anarchy (PoA) is mathematically defined as follows:
\begin{align} 
\label{Eq: PoA}
    PoA &= \frac{C(f_{UE})}{C(f_{SO})} \nonumber \\
    & = \frac{\sum_{a \in A} [ \alpha \cdot t_{a,RN}(x_a^{UE}) ] \cdot x_a^{UE} +  \sum_{m \in \Psi} [ \alpha (WT_{m,a} (q^{UE}) + ST_{m,a}) + C_{m,a} ] \cdot x_{a,m}^{UE}  } {\sum_{a \in A} [ \alpha \cdot t_{a,RN}(x_a^{SO}) ] \cdot x_a^{SO} +  \sum_{m \in \Psi} [ \alpha (WT_{m,a} (q^{SO}) + ST_{m,a}) + C_{m,a} ) ] \cdot x_{a,m}^{SO}  }
\end{align}

{Note that capacitated services, particularly on freeway routes, play a critical role in influencing UE and SO conditions within multimodal transportation networks. When freeway capacities are approached or exceeded, the resulting congestion significantly elevates the costs associated with using these routes, including increased travel times. These capacity constraints can also lead to inefficiencies, particularly in high-demand scenarios, where the costs on freeway routes can far outweigh those of alternative modes. As highlighted in \cite{correa2004selfish}, such conditions create substantial penalties that naturally push commuters towards alternative intermodal options, such as PT or SMSs. This shift helps redistribute demand, reduces pressure on congested freeways, and enhances overall system efficiency. By incorporating multimodal paths, the system can better manage congestion and capacity constraints, thus optimizing the overall flow of traffic.}

Providing an upper bound for the Price of Anarchy is necessary in analytical frameworks due to its capacity to provide insights about the potential negative consequences of not implementing proper traffic management. This upper bound allows policymakers and transportation planners to understand the range of possibilities and make informed decisions about traffic management, infrastructure investment, and policy interventions to strike the right balance between individual freedom and system efficiency in transportation networks. In this matter, the literature commonly discusses upper bounds for two versions of congestion games: atomic and non-atomic congestion games~\citep{bilo_price_2019}. In the atomic scenario~\citep{rosenthal_class_1973}, traffic is discrete, with each player contributing to a congestion level of one. However, in the non-atomic case~\citep{beckmann1956studies},~\citep{pigou2017economics}, each player controls a negligible amount of traffic, resulting in a continuous flow.

We provide a coarse approximation for the upper bound of PoA in multi-modal transportation systems. This is due to the complexity of the link cost function we define. This function is not linear and is bi-variate (link flow variable for travel time and path flow variable for waiting time). Additionally, we are considering the atomic case in this study. All these characteristics make the theoretical analysis of PoA very complex since the application of known upper bounds in the literature or the formulation of a tighter one is not straightforward. However, interested readers may refer to~\cite{bilo_price_2019} for affine (i.e., linear) congestion games,~\cite{qiao_price_2023} for PoA with exponential cost functions, and~\cite{roughgarden_how_2002} for a mathematical and methodological approach to defining bounds for PoA. 

\paragraph{\textbf{Proposition 2}} The price of Anarchy in such a multi-modal system has the following upper bound: 
\begin{equation}
    PoA \leq  B(\beta) + |A| \cdot Q^2 + |A||\Psi | \cdot Q 
\end{equation}

\subparagraph{\textbf{Proof.}} See \ref{app::UB_PoA}.

\section{Numerical Results}
\label{sec::ER}
In this section, we present the numerical results obtained from our study, which aims to analyze the commuter's behavior under different scenarios. We first evaluate the proposed framework on a synthetic network as a proof-of-concept and perform an in-depth analysis to validate the model. Second, we consider a larger multi-modal test case, the Sioux Falls network~\citep{yin_simulation-based_2022}, to validate and discuss the computational performance of the model. By conducting this analysis, we also aim to offer valuable insights into commuters' decision-making process when selecting their mode of transportation and the factors influencing their choices. We solve the UE and SO models using the Gurobi optimizer~\citep{gurobi}. \newref{Table}{tab:2} displays the experiment settings, explicitly mentioning the differing values for the two test cases. It is worth mentioning that for our experiments, we applied both the continuous and integer settings and the results show the same flow pattern for the two settings. Thus, we only report in this section the integer results.

\subsection{Validation and Analysis on a Synthetic Network}
\label{sec:ER1}
The proposed model is implemented for the network illustrated by \newref{Figure}{figNet1}. The network has four nodes, serving as both origins and destinations. Nodes 3 and 4 are transfer nodes to the metro system. We assume a demand of 100 commuters for every OD pair.  

\begin{figure}[!t]
 \centering
 \includegraphics[width=0.8\textwidth]{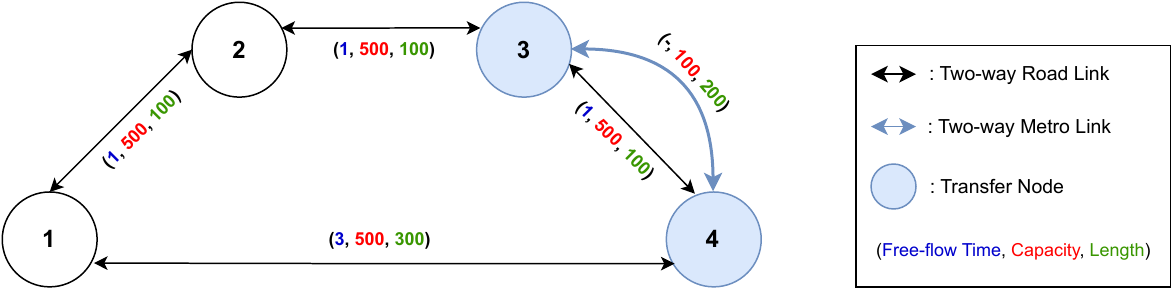}
 \caption{A multi-modal network with four nodes and twelve OD pairs.}
 \label{figNet1}
\end{figure}

We consider three scenarios, illustrated by \newref{Figure}{figTest_Sc}, and analyze the commuters' behavior under UE and SO principles. In Scenario (1), we only consider personal means of transport (car, walking, and biking) and PT (bus and metro). In Scenario (2), we allow SMSs as door-to-door services (without intermodality). Scenario (3) includes all possible travel options, i.e., SMSs and inter-modes are available. For this scenario, the model's solving time is 12 seconds.

\begin{figure}[h!]
 \centering
 \includegraphics[width=0.9\textwidth]{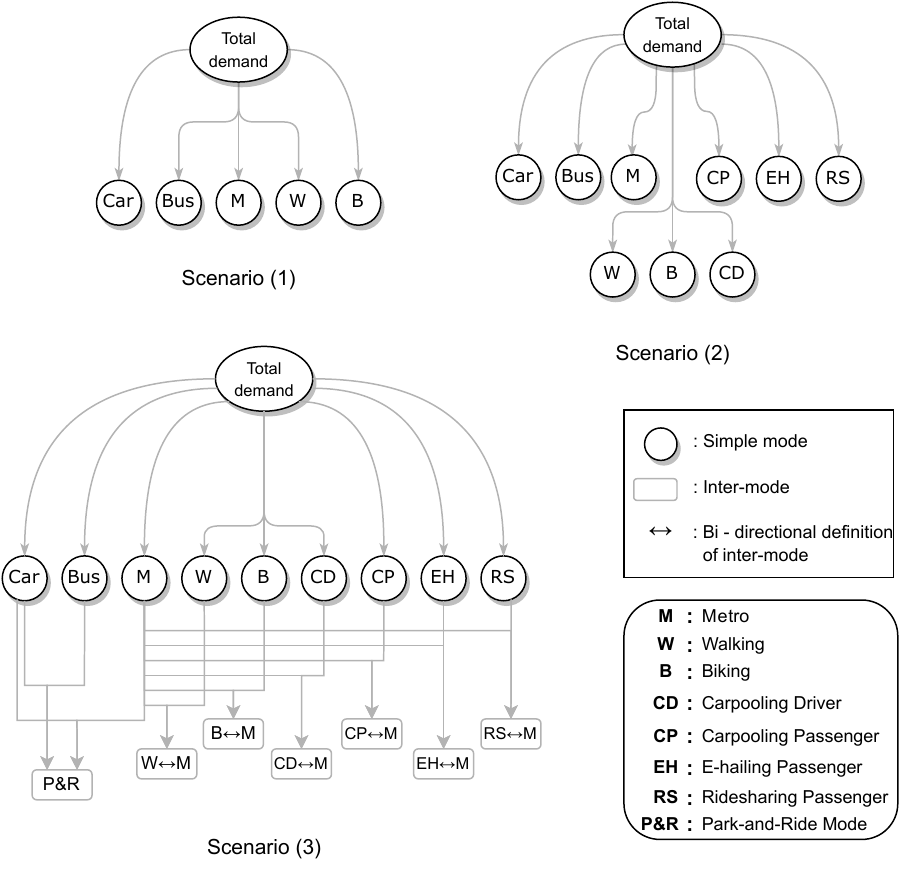}
 \caption{Available travel options for three scenarios.}
 \label{figTest_Sc}
\end{figure}

\newref{Figure}{figTest_Sc2_1} presents the modal share for the first and second scenarios, i.e., the proportion of commuters using each travel mode. Additionally, \newref{Table}{tabSC_1} presents the traffic assignment regarding mode and path choices for the UE solution. {We report the full flow pattern for the first scenario to show that commuters' choices are at equilibrium. However, for the sake of clarity, we only report the chosen options in the remaining scenarios. The non-reported options correspond to unused ones.}

In Scenario (1), $66.67\%$ of commuters choose PT. The bus mode is more frequently used (with $50\%$) due to its accessibility. Since intermodality is not allowed, only commuters between OD pairs $(3,4)$ and $(4,3)$ can access the metro. This demonstrates that using PT when accessible is always preferable to personal means under both UE and SO principles. Let us take the bus service as an example. Bus units operate in the system independently of their passenger flow, contributing simultaneously to traffic congestion in the road network. Besides, a bus can hold more passengers than a car, and since we consider BPR as our travel time function, having as many commuters as possible traveling by PT is the best option. However, it is essential to acknowledge that PT has inherent limitations, such as limited capacity and frequency, restricting its availability to all commuters. Consequently, some commuters will use private cars when PT is not viable. Consider the OD pair \textit{(2,3)} as an example. This pair has no metro links between the origin and the destination. However, there are bus connections between node 2 and node 3. The commuters cannot use the bus mode because it has reached its capacity due to onboard passengers between the OD pair \textit{(2,4)}. 

Moreover, since the network is not congested, solutions under UE and SO are identical, following the observations in~\cite{prashker_observations_2000}. None of the commuters in this scenario have opted for walking or biking due to the higher travel time associated with these modes than other options.

\begin{figure}[!t]
 \centering
 \includegraphics[width=1\textwidth]{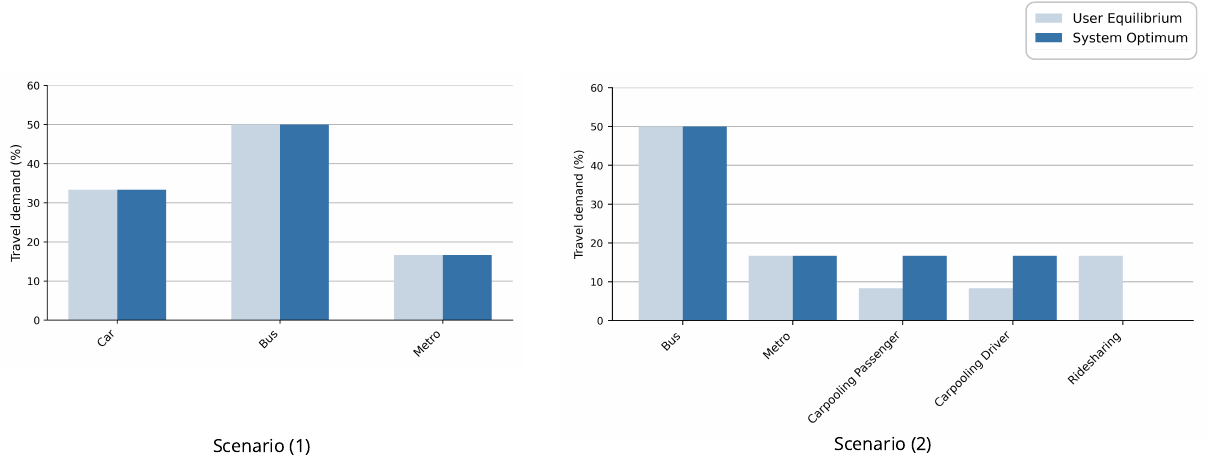}
 \caption{Use of modes in the synthetic test case for the first and second scenarios.}
 \label{figTest_Sc2_1}
\end{figure}

\begin{table}[!t]
\centering
\caption{\textbf{Traffic distribution for the first and second scenarios.}}
\resizebox{\textwidth}{!}{
\begin{tabular}{ c|cccc||c|cccc}
\hline
\textbf{OD pair} & \textbf{Mode} & \textbf{Path} & \textbf{Flow} & \textbf{Gen. Cost*} & \textbf{OD pair} & \textbf{Mode} & \textbf{Path} & \textbf{Flow} & \textbf{Gen. Cost*} \\
\hline
\hline
\multicolumn{10}{c}{Scenario 1: Baseline Network} \\
\hline
\hline
(1,2) & Car & 1-2 & 0 & 12.2 & (1,2) & Car & 1-4-3-2 & 0 & 58.5 \\
(1,2) & Bus & 1-2 & \textbf{100} & 1.4 & (1,2) & Bus & 1-4-3-2 & 0 & 6.4 \\ 
(1,3) & car & 1-2-3 & \textbf{100} & 24.4 & (1,3) & car & 1-4-3 & 0 & 46.3 \\
(1,3) & Bus & 1-2-3 & 0 & 2.9 & (1,3) & Bus & 1-4-3 & 0 & 4.8 \\
(1,4) & car & 1-4 & 0 & 34.1 & (1,4) & car & 1-2-3-4 & 0 & 36.6 \\ 
(1,4) & Bus & 1-4 & \textbf{100} & 3.4 & (1,4) & Bus & 1-2-3-4 & 0 & 4.4 \\ 
(2,1) & car & 2-1 & 0 & 12.2 & (2,1) & car & 2-3-4-1 & 0 & 58.5 \\
(2,1) & Bus & 2-1 & \textbf{100} & 1.4 & (2,1) & Bus & 2-3-4-1 & 0 & 6.4 \\
(2,3) & car & 2-3 & \textbf{100} & 12.2 & (2,3) & car & 2-1-4-3 & 0 & 58.5 \\
(2,3) & Bus & 2-3 & 0 & 1.5 & (2,3) & Bus & 2-1-4-3 & 0 & 6.3 \\
(2,4) & car & 2-3-4 & 0 & 24.4 & (2,4) & car & 2-1-4 & 0 & 46.4 \\
(2,4) & Bus & 2-3-4 & \textbf{100} & 2.9 & (2,4) & Bus & 2-1-4 & 0 & 4.9 \\
(3,1) & car & 3-2-1 & \textbf{100} & 24.4 & (3,1) & car & 3-4-1 & 0 & 46.3 \\
(3,1) & Bus & 3-2-1 & 0 & 2.9 & (3,1) & Bus & 3-4-1 & 0 & 4.8 \\ 
(3,2) & car & 3-2 & \textbf{100} & 12.2 & (3,2) & car & 3-4-1-2 & 0 & 58.5 \\
(3,2) & Bus & 3-2 & 0 & 1.5 & (3,2) & Bus & 3-4-1-2 & 0 & 6.3 \\
(3,4) & car & 3-4 & 0 & 12.1 & (3,4) & car & 3-2-1-4 & 0 & 58.6 \\
(3,4) & Bus & 3-4 & 0 & 1.4 & (3,4) & Bus & 3-2-1-4 & 0 & 6.4 \\
(3,4) & M & 3-4 & \textbf{100} & 3.7 & (4,1) & car & 4-1 & 0 & 34.1 \\
(4,1) & car & 4-3-2-1 & 0 & 36.6 & (4,1) & Bus & 4-1 & \textbf{100} & 3.4 \\
(4,1) & Bus & 4-3-2-1 & 0 & 4.4 & (4,2) & car & 4-3-2 & 0 & 24.4 \\
(4,2) & car & 4-1-2 & 0 & 46.3 & (4,2) & Bus & 4-3-2 & \textbf{100} & 2.9 \\
(4,2) & Bus & 4-1-2 & 0 & 4.9 & (4,3) & car & 4-3 & 0 & 12.1 \\
(4,3) & car & 4-1-2-3 & 0 & 58.6 & (4,3) & Bus & 4-3 & 0 & 1.4 \\
(4,3) & Bus & 4-1-2-3 & 0 & 6.4 & (4,3) & M & 4-3 & \textbf{100} & 3.7 \\
\hline
\hline
\multicolumn{10}{c}{Scenario 2: Multi-modal Network with SMSs} \\
\hline
\hline
(1,2) & Bus & 1-2 & 100 & 1.5 & (1,3) & CP & 1-2-3 & 50 & 5.5 \\
(1,3) & CD & 1-2-3 & 50 & 10.4 & (1,4) & Bus & 1-4 & 100 & 3.4 \\
(2,1) & Bus & 2-1 & 100 & 1.5 & (2,3) & RS & 2-3 & 100 & 8.0 \\
(2,4) & Bus & 2-3-4 & 100 & 3.0 & (3,1) & CP & 3-2-1 & 50 & 5.5 \\
(3,1) & CD & 3-2-1 & 50 & 10.4 & (3,2) & RS & 3-2 & 100 & 8.0 \\ 
(3,4) & M & 3-4 & 100 & 3.7 & (4,1) & Bus & 4-1 & 100 & 3.5 \\
(4,2) & Bus & 4-3-2 & 100 & 2.9 & (4,3) & M & 4-3 & 100 & 3.7 \\
\hline
\end{tabular} \label{tabSC_1} 
}\\
\footnotesize{*: Generalized cost for using the corresponding mode and path between OD pairs.}
\end{table}

In Scenario (2), commuters can choose all modes without intermodality. Like the previous scenario, the bus mode is used to its full capacity. The metro is also preferred when accessible. For SMSs, $16.66\%$ choose carpooling (both as passengers and drivers), and $16.67\%$ use ridesharing. \newref{Table}{tabSC_1} shows that ridesharing is used for short-distance trips (OD pairs \textit{(2,3)} and \textit{(3,2)}) while carpooling is chosen for long-distance trips. This can be explained by the fact that ridesharing has a high trip fare, increasing with the length of the trip (traversed links). In this scenario, no commuter is willing to use the car or the e-hailing service because, without considering personal preferences, it is always beneficial to share the cost of the trip with a passenger through carpooling or ridesharing.

Note that under the UE conditions and for the same OD pair, the cost of carpooling passengers is not the same as that of drivers. This demonstrates the redefinition of the UE conditions under capacity constraints. For example, consider the OD pair \textit{(1,3)}. For these commuters, the bus is the option with minimum cost. Nevertheless, it reaches full capacity because of commuters between \textit{(1,2)}. The second best option is carpooling as a passenger. However, this is feasible only if there are enough drivers, so only half of these commuters can choose the CP mode. The cost for a carpooling driver is high due to the cost of using the car, but since this cost is less than driving alone or ridesharing, this option is used by the other half of commuters.

When traffic is assigned under SO, the ridesharing service is not used. Our SO objective function does not consider the ridesharing service provider's profits. Thus, from a system point of view, the goal is to reduce both the passengers' and drivers' costs. If the number of CP increases, the number of CD also increases accordingly, which will attract the demand to switch from the ridesharing service to the carpooling service.

\newref{Figure}{figTest_Sc2_2} presents the use of the modes for the third scenario, where intermodality is available. In this case, we observed that, with the same model configuration as the previous two scenarios, more than 80\% of commuters will use PT due to its accessibility, while the rest of the commuters will use carpooling. However, to present a more realistic and complete analysis, we reduce the bus frequency ($freq_{Bus} = 2$) and provide the traffic assignment information in \newref{Table}{tabSC_2}. In this case, when traffic is assigned based on UE, $36.67$\% of commuters use PT, 10.83\% carpool as passengers for door-to-door trips and not for intermodality, and 10.08\% participate in ridesharing, while $15.75\%$ of commuters use the CD mode. These drivers benefit from the possibility of having a passenger on board (CP or I(CP)) for part of or the entire trip to compensate for the high cost of using their car. Intermodality is used with SMSs, representing $26.66\%$ of the trips (neither P\&R mode nor walking and biking for the first/last mile). Similar to the previous scenario, ridesharing is mainly used for short trips (either as a door-to-door service or combined with the metro for a short distance). In this case, ridesharing is used instead of carpooling due to the high demand for the carpooling service and, thus, a high waiting time. The trip with intermodality at the UE showed that our model could capture these complex trips.  

Under the SO principle, the ridesharing demand is directed toward carpooling. Here, many commuters use M\&CP because using PT improves the system's overall objectives. This increase in the M\&CP mode attracts more drivers to participate in carpooling (CD mode) to reduce their costs since they can transport passengers from the transfer node to their destination.

\begin{figure}[!t]
 \centering
 \includegraphics[width=1\textwidth]{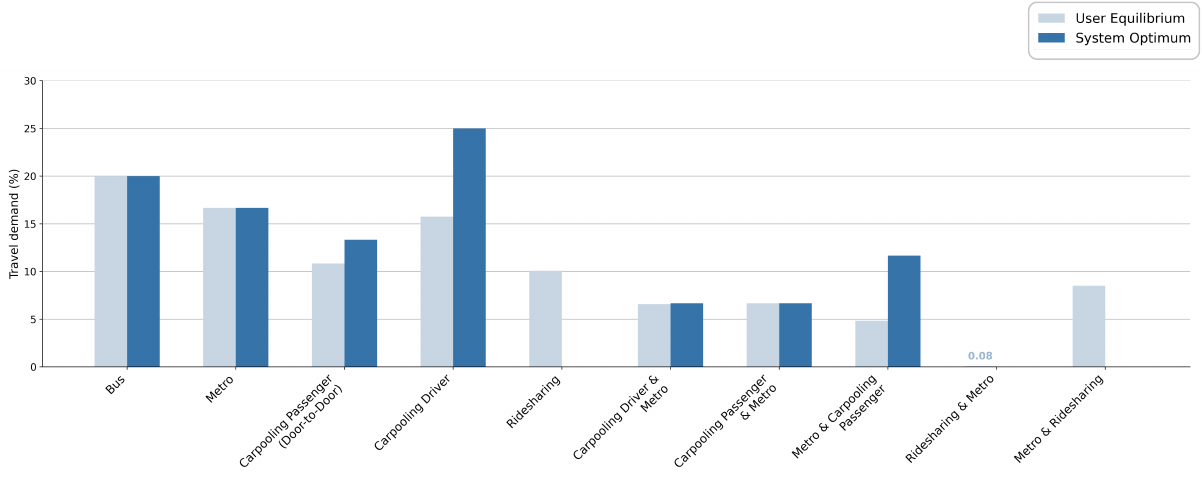}
 \caption{Use of modes in the synthetic test case for the third scenario. We distinguish the use of SMSs for the first mile or last mile (e.g., Ridesharing \& Metro refers to the use of ridesharing for the first mile, while Metro \& Ridesharing refers to the use of ridesharing for the last mile of the trip)  }
 \label{figTest_Sc2_2}
\end{figure}

\begin{table}[!t]
\centering
\caption{\textbf{Traffic distribution for the third scenario.}}
\resizebox{\textwidth}{!}{
\begin{tabular}{ c|cccc||c|cccc}
\hline
\textbf{OD pair} & \textbf{Mode} & \textbf{Path} & \textbf{Flow} & \textbf{Gen. Cost} & \textbf{OD pair} & \textbf{Mode} & \textbf{Path} & \textbf{Flow} & \textbf{Gen. Cost} \\
\hline
\hline
(1,2) & Bus & 1-2 & 40 & 1.6 & (1,2) & RS & 1-2 & 60 & 11.0 \\
(1,3) & CP & 1-2-3 & 80 & 7.3 & (1,3) & CD & 1-2-3 & 20 & 10.4 \\
(1,4) & Bus & 1-4 & 40 & 3.5 & (1,4) & CD\&M & 1-2-3-4 & 60 & 14.2 \\
(2,1) & Bus & 2-1 & 40 & 1.6 & (2,1) & RS & 2-1 & 60 & 11.0 \\
(2,3) & Bus & 2-3 & 40 & 1.6 & (2,3) & CD & 2-3 & 60 & 5.2 \\
(2,4) & CD\&M & 2-3-4 & 20 & 8.9 & (2,4) & CP\&M & 2-3-4 & 80 & 7.4 \\
(3,1) & CD& 3-2-1 & 100 & 10.4 & (3,2) & Bus & 3-2 & 40 & 1.6 \\
(3,2) & CP & 3-2 & 30 & 3.6 & (3,2) & CD & 3-2 & 30 & 5.2 \\
(3,4) & M & 3-4 & 100 & 3.7 & (4,1) & Bus & 4-1 & 40 & 3.5 \\
(4,1) & M\&CP & 4-1 & 60 & 11.0 & (4,2) & M\&RS & 4-3-2 & 100 & 14.8 \\
(4,3) & M & 4-3 & 100 & 3.7 & & & & &  \\
\hline
\end{tabular} \label{tabSC_2} 
}
\end{table}

\newref{Figure}{figNet1_PA}(A) presents the total cost evolution under both UE and SO, as well as the price of anarchy (PoA) for the road network, defined by Equation (\ref{Eq: PoA}) when the total demand increases. The total cost under UE is approximately 1.2 times higher than SO when the demand is increased by a factor of 10. This means that when the network is congested, commuters' personal decisions decrease the overall performance level of the transport system. It is expected due to the UE objective function. In such a case, increasing the capacity of PT or promoting more carpooling through subsidies, for example, may help bring the solution under UE to that under SO. Additionally, from these numerical experiments, we observed that the PoA is still bounded by the upper bound for general atomic congestion games, with linear link cost functions, shown to be equal to $2.5$~\citep{christodoulou_price_2005}. 

\begin{figure}[!h]
\centering
 \includegraphics[width=0.8\textwidth]{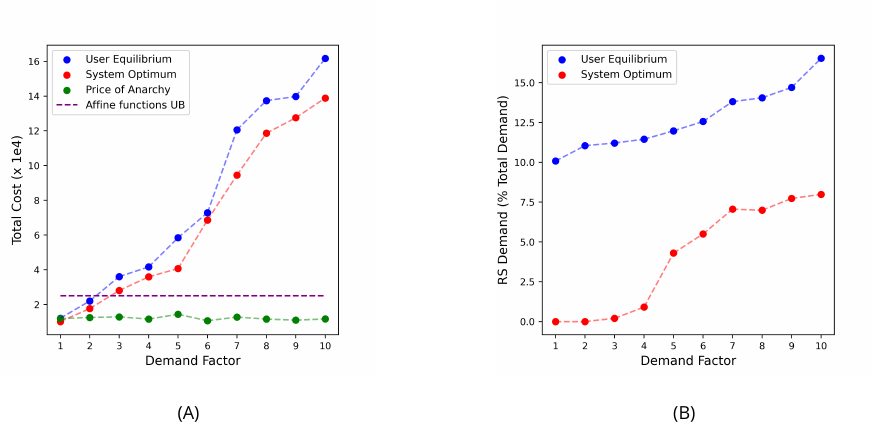}
\caption{\textbf{(A)}: Price of Anarchy - Evolution of the total cost in the synthetic test case; \textbf{(B)}: Evolution of ridesharing use for the synthetic test case, w.r.t the total demand.}
\label{figNet1_PA}
\end{figure}

Moreover, we investigated the variation in ridesharing use when the total demand increases. \newref{Figure}{figNet1_PA}(B) illustrates, under UE, the increased use of ridesharing services, which can be explained by the fact that, for carpooling drivers, ridesharing becomes more attractive since the travel time is high and the difference in terms of monetary cost between using their car and ridesharing becomes comparable. Besides, under SO, since the PT reaches its total capacity, demand is directed to carpooling and ridesharing. However, the ridesharing service is more present due to the high waiting time for carpooling passengers. To further investigate the use of carpooling and ridesharing services, we analyze the impact of distance on passengers' choices. 

\subsubsection{Carpooling versus Ridesharing Service}
The previous results show that the carpooling service is preferred over ridesharing due to its cost-effectiveness for both drivers and passengers. In particular, drivers reduce their costs by participating in carpooling as long as a passenger is available (for instance, carpooling drivers between OD pair $(2,3)$ having passengers on board between OD pair $(2,4)$ for part of their trip). However, studies such as those by~\cite{yan_stochastic_2019},~\cite{gheorghiu_for_2018} and~\cite{vanoutrive_what_2012} have shown that drivers are not attracted by the carpooling service for short-distance trips. In line with these works, we modified the proposed model by restricting the carpooling service to be available only for trips longer than a pre-defined threshold for short distances. For the synthetic network, we set this value to $100$ distance units, corresponding to the minimum link length in this test case, and analyzed the variation in the use of modes with this configuration. 

\newref{Figure}{figNet1_CPres} illustrates the use of modes when short-distance carpooling is not allowed. The results show that the modal share has changed compared to the third scenario previously described, with a decrease in the carpooling use and an increase in the ridesharing percentage. There are fewer carpooling drivers (either door-to-door, or using intermodality), under both UE and SO principles. Fewer carpooling drivers in the system also means fewer carpooling passengers. This led to an increase in the ridesharing service since all the demand that was previously fulfilled by carpooling for short distances, as well as short-distance carpooling to access or egress a PT station (intermodality), was directed towards the ridesharing service due to its effectiveness and the saturation of the bus network. It should be noted that the main purpose of this restriction is to consider realistic scenarios, and its impact on commuters' behavior needs to be further analyzed on larger networks.

\begin{figure}[!h]
\centering
\includegraphics[width=1\textwidth]{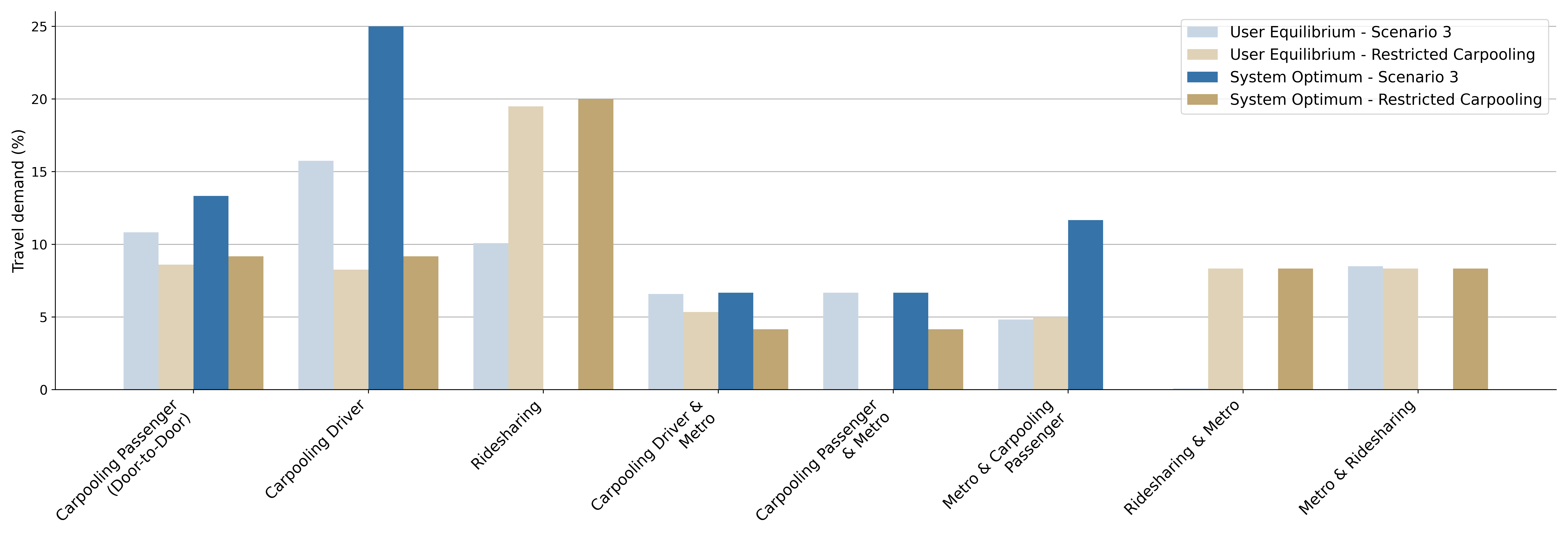}
\caption{Use of modes in the synthetic test case when short-distance carpooling is not allowed}
\label{figNet1_CPres}
\end{figure}

\subsubsection{Sensitivity Analysis}
In this section, we investigate the sensitivity of our model to key input parameters that significantly influence travel behavior. Our objective is to gain a deeper understanding of factors that can impact travel choices, focusing on the pricing of different travel modes and the concept of the value of time. By examining the sensitivity of our model to changes in these parameters, we shed light on the robustness of our analytical framework.

\paragraph{\textbf{Trip Fares for Shared Mobility Services}} 
The monetary cost of a travel option is an important factor in the decision-making process. We tested the sensitivity of our model to this factor in terms of the use of the modes in three different cases. The first case is the baseline (previously presented), in which the carpooling service is less expensive than ridesharing and car modes. The second case is a scenario in which ridesharing is less expensive than other travel options ($TF_{RS} = 0.9 \; ; \; TF_{CP} = 1.2$), while in the last case, the car is the most attractive travel option ($\beta = 0.1 \; ; \; PF_{car} = 0.5$). \newref{Figure}{figSA_1}(A) illustrates the use of the modes when traffic is assigned based on UE, and \newref{Figure}{figSA_1}(B) is for the SO principle. Under both UE and SO principles, most of the travel demand for carpooling is directed towards ridesharing services when the latter is less expensive. The remaining carpooling demand is due to the high waiting time for the ridesharing service. In this case, the ridesharing service mostly competes with PT since it attracts door-to-door requests. This result is in line with the work of~\cite{zhu_analysis_2020}. On the other hand, when the car is more attractive, all SMS passengers prefer to use their own car to optimize their travel costs, as expected. However, from the system point of view (SO principle), the optimal strategy for this scenario is to have more drivers willing to participate in carpooling, which reduces the number of cars in the network. Thus, all commuters will experience shorter travel times (due to less congestion). 

\begin{figure}[h!]
 \centering
 \includegraphics[width=1\textwidth]{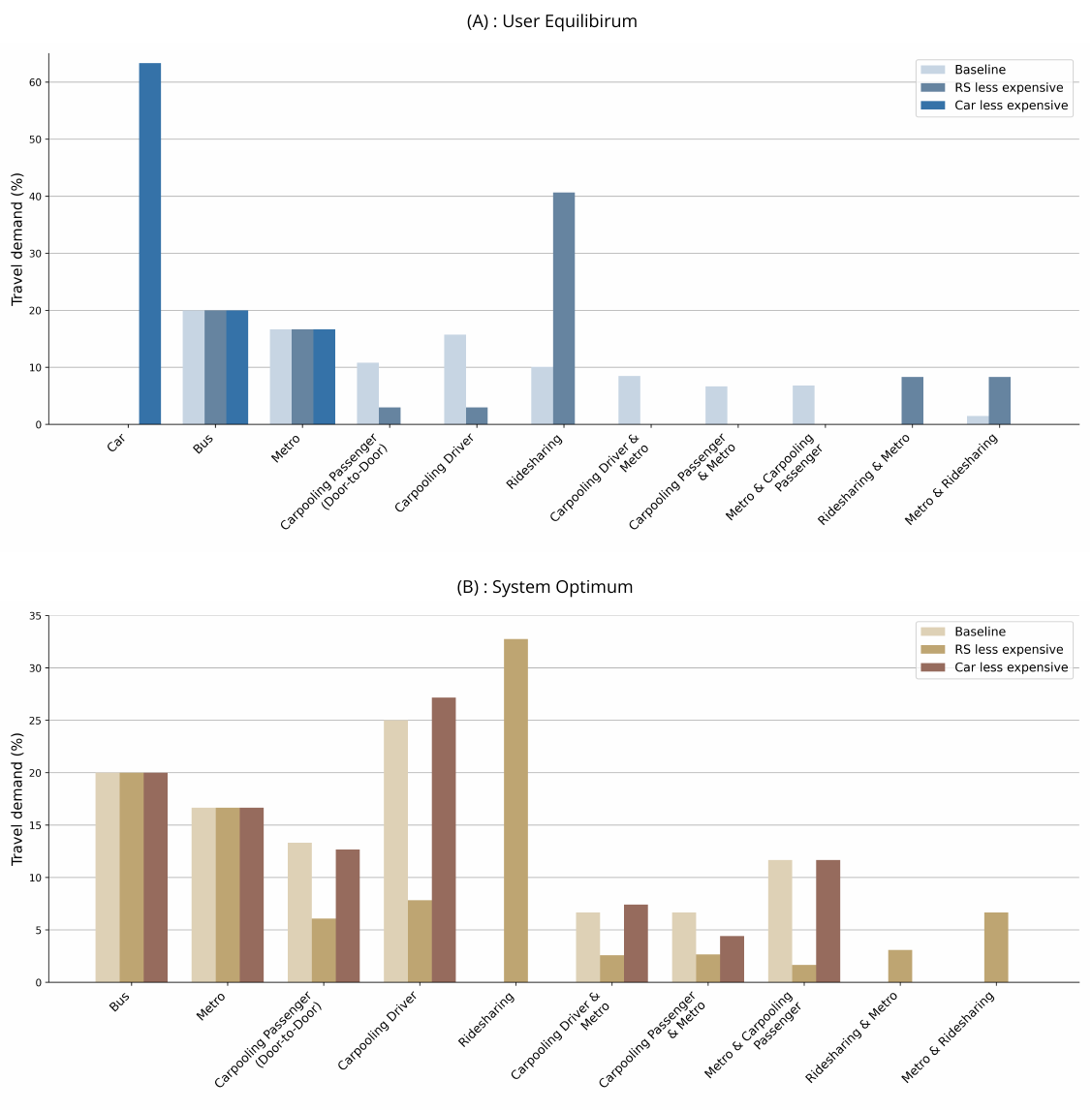}
 \caption{Sensitivity of use of modes to trip fares for shared mobility services.}
 \label{figSA_1}
\end{figure}

\paragraph{\textbf{Trip Fare for Public Transport}}
Public transport is an attractive travel option with a low monetary cost and a constant waiting time. In what follows, we increased the trip fare for the bus mode and analyzed the changes in the travel behavior. We focus only on the bus mode because of its interactions with the car mode and shared mobility services. We consider five different trip fares for the bus ( $TF_{Bus} \in \{0.3 \;,\; 0.7 \;,\; 1 \;,\; 5 \;,\; 7\}$ ), and show the use of the modes in \newref{Figure}{figSA_2}. When the increase in the trip fare is not significantly high, commuters do not change their travel behavior. This happens because, in our model, the PT operates under a fixed schedule. In other words, if the initial bus passengers switch to an SMS due to a higher price, for instance, they will increase their travel time due to the congestion produced by their flow, as well as the flow of empty buses. Thus, even if the bus trip fare increases, they will keep using this mode to avoid creating congestion and having a higher travel time. However, when the compromise between the travel time and monetary cost starts to be significant (i.e., when the trip fare of buses is fairly expensive), the travel demand for the bus mode is mostly directed towards door-to-door SMSs. 

\begin{figure}[h!]
 \centering
 \includegraphics[width=1\textwidth]{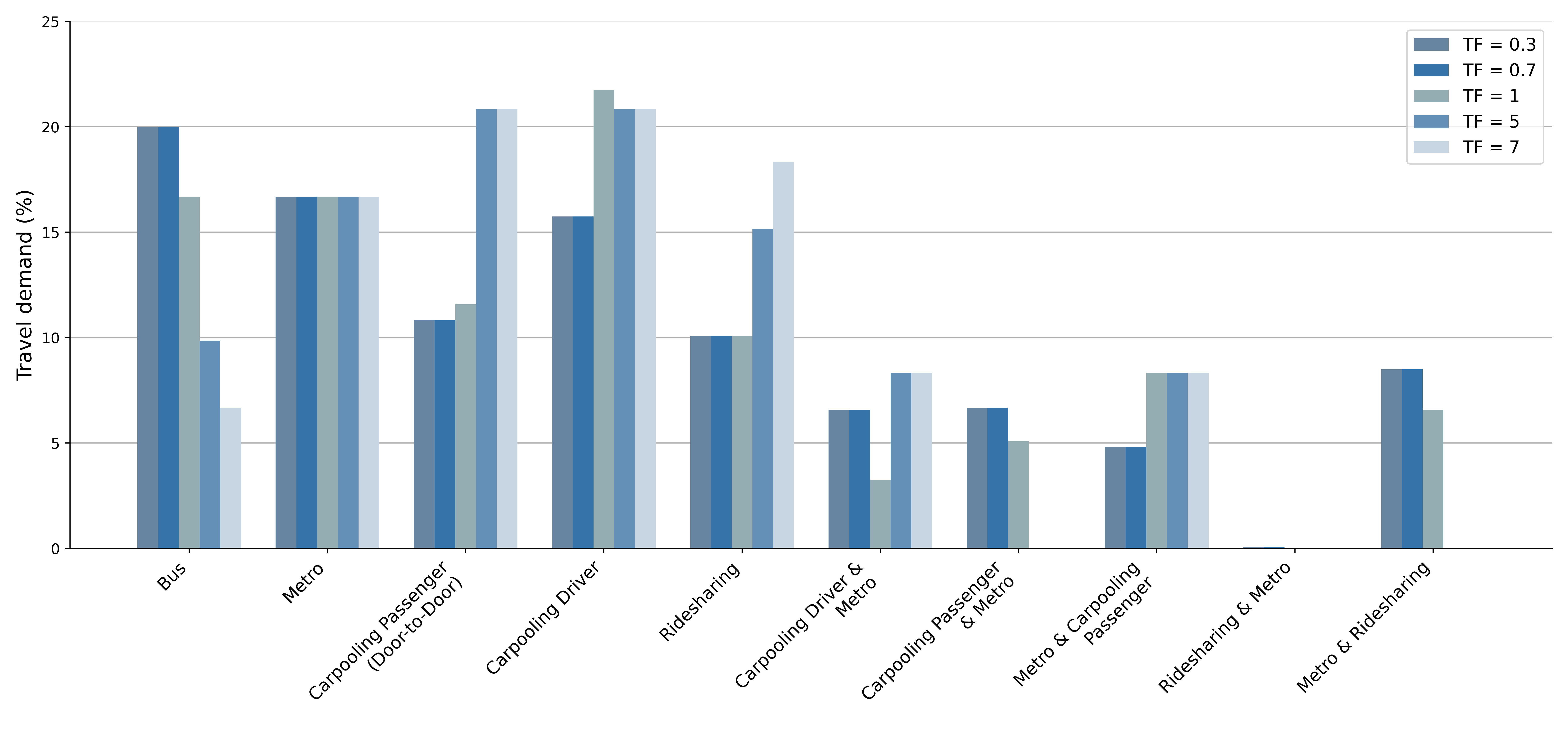}
 \caption{Sensitivity of use of modes to the trip fare for public transport.}
 \label{figSA_2}
\end{figure}

\paragraph{\textbf{Value of Time}}
The value of time (VOT) determines the importance of time for commuters. The higher this value is, the more expensive a unit of time is. In \newref{Figure}{figSA_4}, we present the use of modes for three different values of time. On the one hand, when a low VOT is used, commuters are more sensitive to the monetary cost of the travel options. Thus, they are more willing to bike, either for the whole trip or as a first/last mile option, absorbing the demand for SMSs. On the other hand, a high VOT leads to a shift to using private cars and door-to-door SMSs, rather than intermodality, which implies a high waiting time due to the PT transfers. In this case, all passengers using the ridesharing service share the vehicle with commuters having the same OD pair to avoid making a detour to pick up or drop off passengers with a different origin or destination.

\begin{figure}[h!]
 \centering
 \includegraphics[width=1\textwidth]{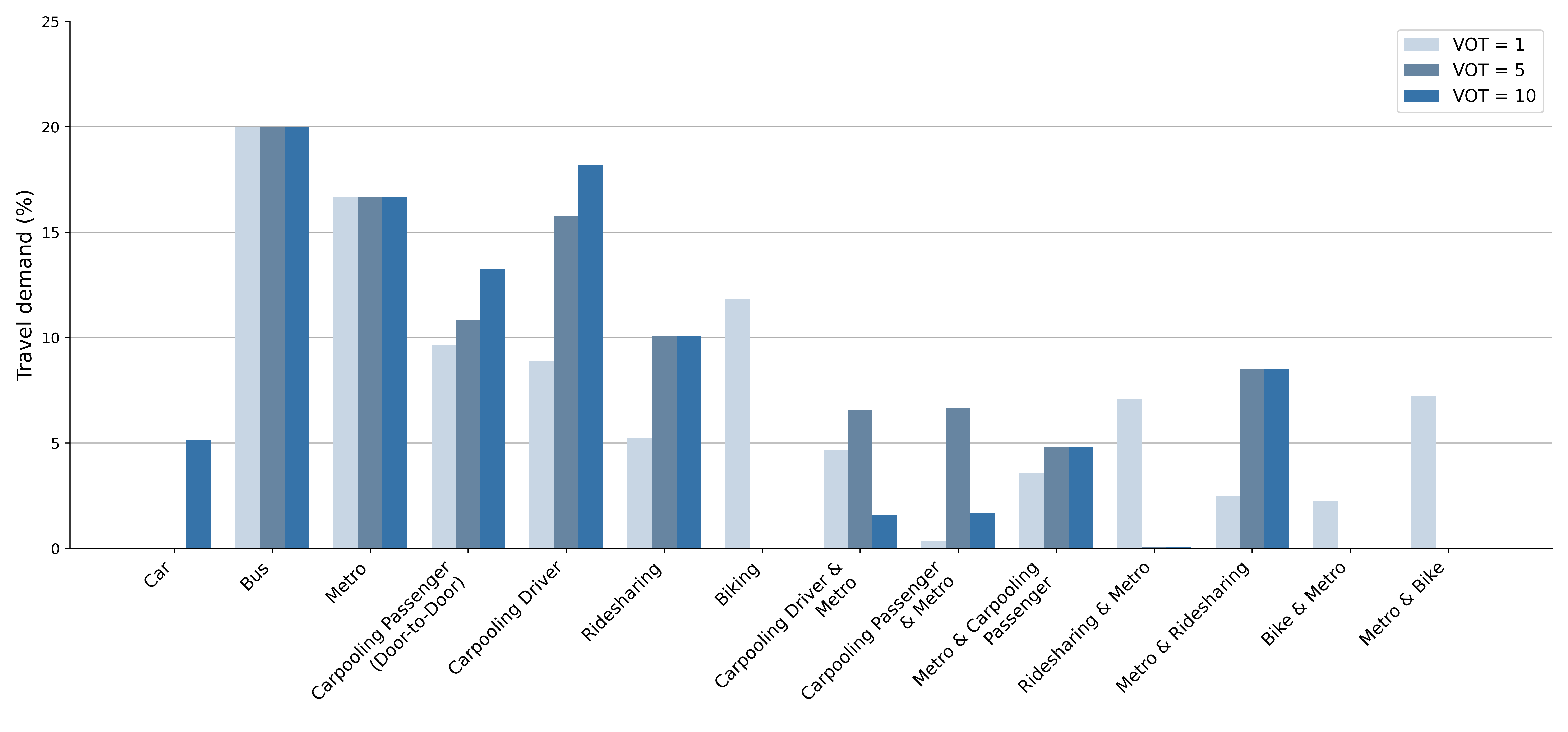}
 \caption{Sensitivity of use of modes to the value of time.}
 \label{figSA_4}
\end{figure}

\subsection{Application and Discussion}
\label{sec:ER2}
In this section, we provide an extensive evaluation and discussion of the applicability of our model to large-scale multi-modal urban transportation networks. We first apply the model to the Sioux Falls network and provide an analysis of the commuters’ travel behavior. Then, we discuss the computational performance of the model, highlighting its resource requirements, scalability, and potential optimization strategies to accommodate larger networks.

\subsubsection{Travel Behavior in Sioux Falls Network}
We analyzed the travel behavior in the multi-modal Sioux Falls network, with PT infrastructure extracted from~\cite{yin_simulation-based_2022} and illustrated by \newref{Figure}{figNet2}. The road network parameters and the travel demand are extracted from~\cite{stabler_transportation_nodate}. The network comprises 24 nodes, with nodes $4$, $6$, $15$, and $19$ as transfer nodes. For the travel demand, we consider the network's 34 most crowded OD pairs \citep{stabler_transportation_nodate}.  

\begin{figure}[h!]
 \centering
 \includegraphics[width=0.6\textwidth]{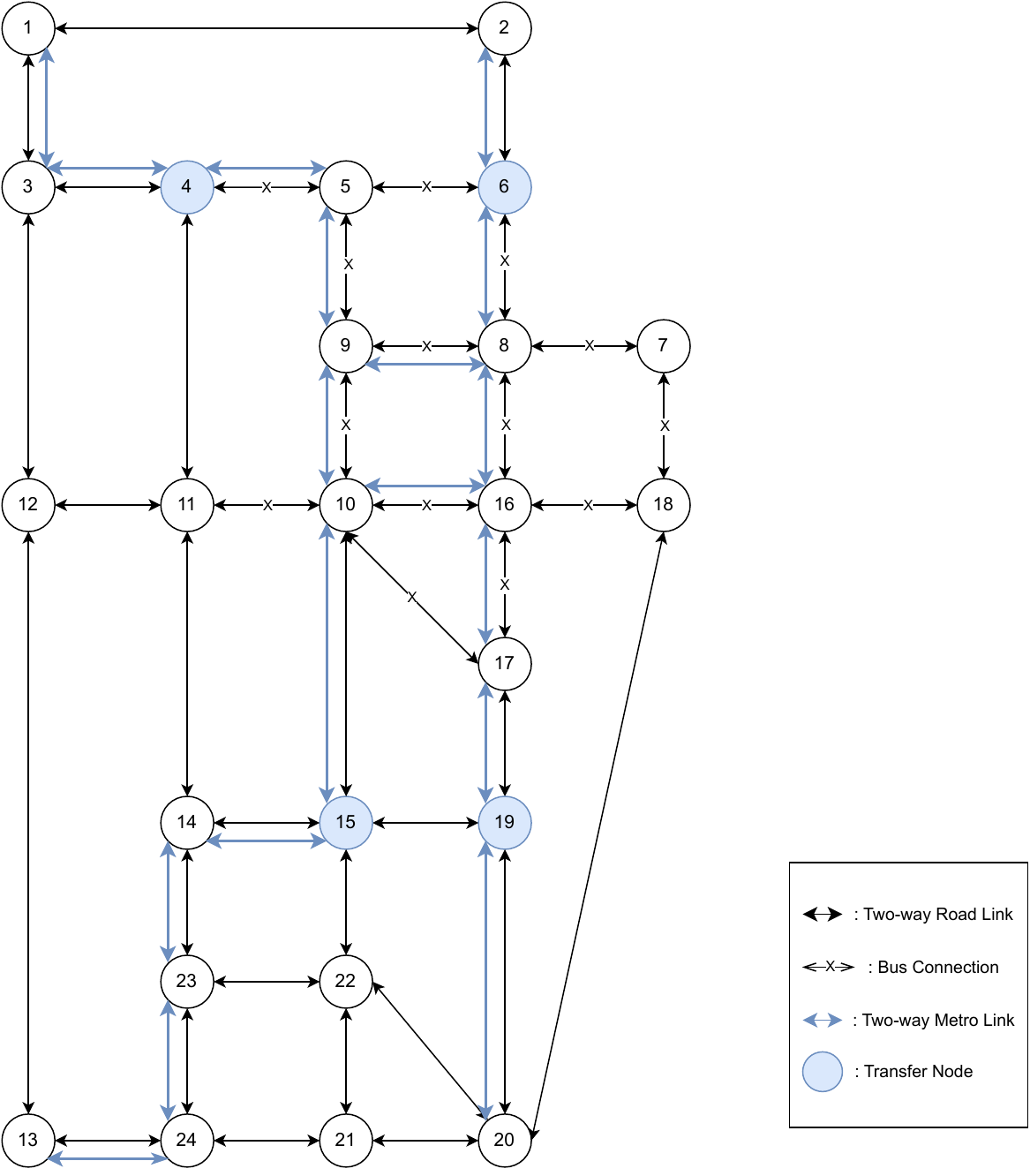}
 \caption{Sioux Falls Network with PT infrastructure.}
 \label{figNet2}
\end{figure}

\newref{Figure}{figSF_sc} illustrates the proportion of commuters using each available mode. We can observe a globally similar travel behavior to that described in \newref{Section}{sec:ER1}. Public transport is widely used under both UE and SO principles. The high mode share for PT is due to its accessibility for the OD pairs considered in this scenario. In particular, it can be seen that most of the pairs are directly connected through PT (either with buses or metro lines), which corresponds to 70\% of the total demand. For OD pairs that are not directly connected by public transport or those whose PT connections are saturated, commuters use shared mobility services. Like the previous case study, carpooling is mostly used for long-distance trips, while ridesharing is used for short-distance trips.

For instance, let us consider the OD pair $(10, 11)$. Since this pair has no direct metro connection, the bus mode is used until its saturation, satisfying 40\% of the demand. Based on the road structure, the link length between node $10$ and node $11$ equals five distance units, which is considered a short distance. Thus, 40\% of commuters use the ridesharing service to reach their destination, and 19\% participate in carpooling through this direct link. However, due to the congestion of this link caused by the traffic between this OD pair, as well as the traffic between nodes $10$ and $12$, the remaining 1\% of the travel demand for $(10,11)$ uses intermodality, taking the metro link between $10$ and $15$, then participating as carpooling passengers between $15$ and $11$.

\begin{figure}[h!]
 \centering
 \includegraphics[width=1\textwidth]{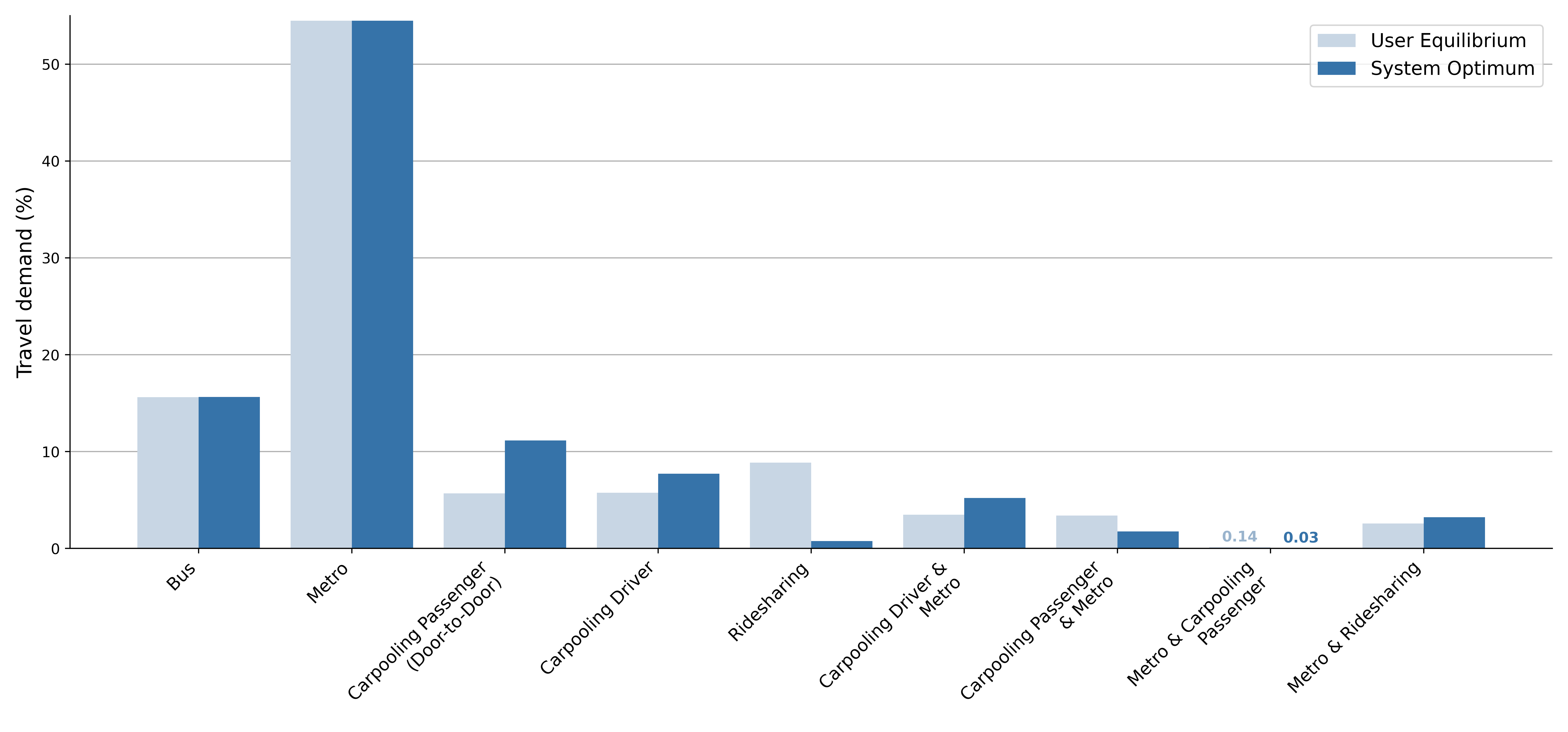}
 \caption{Use of modes in the multi-modal Sioux Falls network, with intermodality.}
 \label{figSF_sc}
\end{figure}

\newref{Figure}{figSF_sc} also shows that under the UE principle, 9.57\% of commuters rely on intermodality with SMSs for the first and last mile, compared to 10.21\% for SO. This increase in intermodality under the SO principle is due to more carpooling drivers. Specifically, when PT is saturated, carpooling is the most attractive option for passengers despite its waiting time, and since the total cost is minimized under the SO principle, the proportion of carpooling drivers will increase to satisfy more carpooling passengers. In other words, even if the cost of using the car is relatively high, more carpooling drivers are willing to accept higher costs under the SO principle to allow other commuters to participate in carpooling as passengers and, thus, decrease the overall system’s cost. However, these carpooling drivers will still try to decrease their own costs by choosing intermodality (CD\&M), consequently reducing the distance on which they use their car.

\newref{Figure}{figSF_congestion} illustrates the impact of intermodality on vehicular (i.e. road) flow across the network, comparing the scenario where intermodal travel is restricted (\newref{Figure}{figSF_congestion}-A) to the scenario where all modes can be combined (\newref{Figure}{figSF_congestion}-B). The results show that allowing intermodality helps reduce congestion, particularly for OD pairs that benefit from PT transfers. For example, congestion around node 23 is eliminated in scenario B, as travelers to and from this node shift to metro and transfer at node 15, using an SMS for the first or last mile. Additionally, intermodality contributes to easing congestion around node 10, the most crowded node representing the city center, highlighting its potential to improve overall network conditions.

\begin{figure}[h!]
 \centering
 \includegraphics[width=0.9\textwidth]{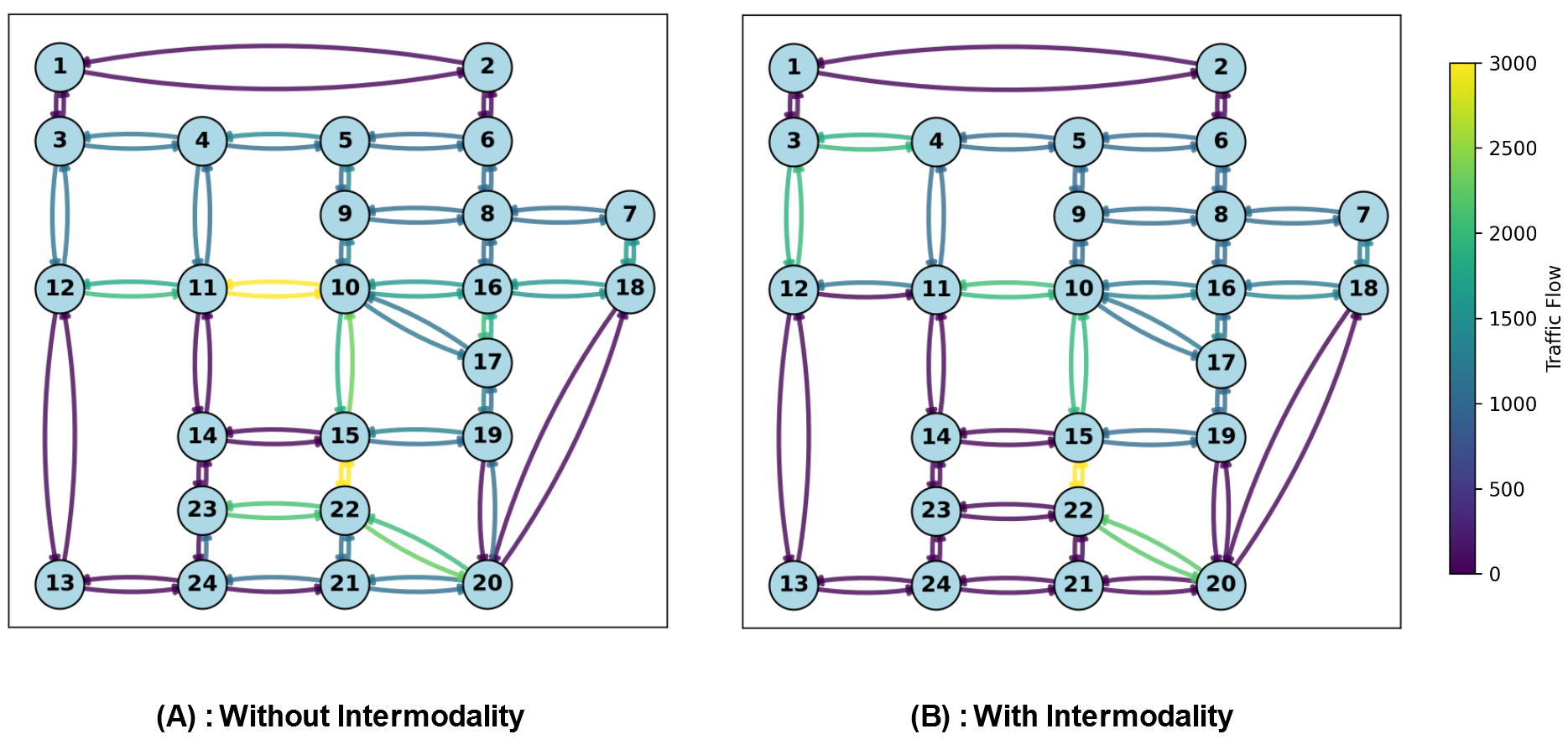}
 \caption{Vehicular link flow distribution in Sioux Falls Network illustrating the impact of intermodality on traffic congestion.}
 \label{figSF_congestion}
\end{figure}

The previously mentioned observations suggest that improving intermodality can improve the total travel cost in such multi-modal systems. In particular, motivating more drivers to participate in carpooling, either for door-to-door trips or through intermodality, will give the transportation planners opportunities to shift from the UE state towards the SO state. 

\subsubsection{Computational Performance}
It is important to note that the computational time for our model exhibits significant variance across different instances. This variability can be attributed to several factors, particularly the size of the transportation network, the public transport structure and its accessibility, the volume of commuters associated with each OD pair, as well as factors affecting the costs of commuters, such as the geographical distance between the nodes and the different input parameters, including the value of time and trip fares. The positioning of the OD pairs within the network is also an important factor since it determines the sharing possibilities for commuters using carpooling or ridesharing services, and these sharing possibilities influence the number of variables of the model. These diverse parameters collectively contribute to the fluctuating computational time experienced during our analysis. Moreover, considering the multimodal and intermodal aspect of the framework, including SMSs, makes the model very challenging to solve on large instances. This is also the case of the relevant works in the literature, such as~\citep{du_modeling_2022} and ~\citep{pi_general_2019}.

Furthermore, the complexity of the solving process in the Gurobi optimizer also influences the performance of our model. Gurobi uses the branch and bound method for solving integer programs. The branching is done on each variable; thereby, the corresponding relaxed program is solved~\citep{junger_mixed_2013}. The complexity of solving a bilinear program is approximated to $O(P^3)$, \textit{P} being the number of decision variables~\citep{nocedal1999numerical} (see~\cite{cimini_exact_2017} for a more detailed and precise study). The complexity of the branch and bound algorithm is $O(Mb^d)$, where M is the cost of expanding sub-problems, b is the branching factor (in our case, the branching is done on each variable, thus $b=2$), and d is the search depth~\citep{morrison_branch-and-bound_2016}. Consequently, for our case study, we can approximate the Gurobi complexity by $O(P^3 2^P)$. However, this is a rough estimation since the existing solvers usually apply many heuristics and internal modifications to optimization problems to improve the solving process.

To estimate and discuss the computational performance of our model, we performed several executions on a subset of 30 randomly chosen OD pairs from a total of 552 for the Sioux Falls network. Experiments were conducted on a 64-bit computer with Intel i7 CPU 2.90GHz and 16GB RAM. We report an average solving time of 4.54 hours with a standard deviation of 6.41 hours. The minimum solving time is 40 seconds for an instance with around 100 commuters for each OD pair, while the maximum time is experienced for the complete scenario, with 19 hours of solving time.  

It is worth mentioning that, even though the reported execution times are high, the model remains suitable for analyzing the use of travel modes in such complex multi-modal transport systems. Essentially, the model is capable of providing, through a single execution, all the necessary elements for this analysis, namely, the optimal traffic assignment for all OD pairs, the use of different available modes, the choice of transfer stations for intermodality with public transport, as well as the optimal matching for passengers and drivers for SMSs. This demonstrates its applicability for day-to-day travel analysis, allowing for advance planning of transport services.

However, we intend to improve this aspect of the model in future work. Some perspectives include solution methods or heuristics suitable for the considered problem since using the existing traffic assignment approaches, such as the Method of Successive Averages (MSA), is inadequate due to the higher complexity of SMSs and the integer nature of traffic flow in our case. Additionally, changing the network's aggregation level is another interesting perspective. In particular, considering the graph nodes as distinct and small regions instead of road intersections (as is currently the case) and thus creating a hypergraph for the real transportation system will eventually allow the application of the model to large-scale urban transportation networks. Finally, using other commercial solvers that might be more efficient than Gurobi for solving integer bilinear problems can also be examined.

\section{Conclusion}
\label{sec::conclusion}
This research proposed a comprehensive mathematical assignment model for multi-modal urban transportation networks, including shared mobility services (SMSs) and public transport (PT). The model is able to address the impacts of SMSs on transportation systems and commuters' choices when integrated with PT. Existing studies have shown shortcomings in providing an integrated framework for evaluating, combining, and optimizing available travel options effectively, leading to potential inefficiencies in planning and higher expenses. In response, our proposed model, where commuters can seamlessly select and combine various travel modes, including SMSs and PT, allows for optimized travel choices.

Two fundamental assignment principles, user equilibrium (UE) and system optimum (SO), are considered for the model formulation (\newref{Program}{P1}). The model not only includes transit, personal, and shared transportation modes but also captures the intermodality in which users can combine more than one travel mode. To the best of our knowledge, this is the most complete assignment model for multi-modal networks (see \newref{Table}{tab:1}). We examined commuters' travel decisions under UE and SO principles and provided an analysis of the price of anarchy in such systems. Our findings reveal carpooling is preferred over ridesharing when PT is not directly accessible. However, in situations of high demand and network congestion, ridesharing becomes more attractive due to longer waiting times for carpooling passengers and reduced costs between cars and ridesharing services because of traffic congestion. Additionally, the carpooling service attracts more commuters for long-distance trips than the ridesharing service. However, the latter becomes more attractive when the fare is low, thus competing with PT by providing door-to-door services to passengers. Furthermore, improving intermodality and encouraging more commuters to participate in carpooling, especially as drivers, will potentially shift the travel behavior from user equilibrium to a system optimum state.

As an initial step toward a unified framework incorporating SMSs and PT, our proposed model offers a baseline for analyzing commuters' behavior within a multi-modal urban transportation network. It facilitates the investigation of planning scenarios and a better understanding of the integrated transportation system. 
{It is important to note that this study primarily focused on modeling traffic assignment within a multimodal transportation system under a specified supply scenario. However, the definition of the optimization variable could be extended to include public transportation \citep{daganzo2019public}, shared mobility services \citep{daganzo2019general}, or a combination of both \citep{liu2021mobility}. Such extensions would inevitably increase the complexity and computational time of the solution significantly. Building on the work of \cite{kim2024computational}, who proposed a framework that integrates travel demand and network assignment models through machine learning methods to enhance consistency and computational efficiency, this approach presents a promising avenue for the further development of urban multimodal systems.}

Other future research directions should consider personal preferences and social criteria in commuters' decision-making processes, the attractiveness of different modes, incorporate service providers' objectives and integrate market equilibrium into the model's formulation. The computational performance of the model also needs to be improved by providing heuristic solutions or different approaches to reduce the size of the network. In addition, including traffic dynamics can enhance the model while increasing the system's complexity and solution calculation. The simulation-based approach can be applied while it makes the system intractable for large-scale applications. Formulating other types of equilibria, such as stochastic user equilibrium, can also be investigated and will require the modification of the objective function to consider the stochastic terms while keeping the same constraints and link costs as the proposed model.



\section{Acknowledgements}

M. Ameli acknowledges support from the French ANR research project SMART-ROUTE (grant number ANR-24-CE22-7264).

\section*{Disclosure of interest}
The authors report there are no competing interests to declare.

\end{sloppypar} 

\typeout{}
\bibliography{references}

\appendix
\section{Linearization of the Objective Function} 
\label{app::linear}
In order to keep the objective function in Equation 1 bilinear, and to be able to use state-of-the-art solvers, we only need to apply the linearization procedure to the integral term with BPR function. The other terms in the objective function are already linear. The linearization procedure is further described in \cite{amini_braess_2022} and \cite{Wei_efficient_2019}.

Let us consider the program T as follows, with $t_{a,n}$ as the BPR function, in which $t_a^0$ corresponds to the free-flow travel time on link $a$, and $c_a$ to the maximum capacity of link $a$. $\eta$ and $\beta$ are the parameters of the BPR function.

\begin{align}
\label{app1_eq1}
    T & =  \min \sum_{a \in A_n}  \int_0^{x_a}  \alpha \cdot t_{a,n}(\omega) \; d\omega = \sum_{a \in A_n} \int_0^{x_a}  \alpha \cdot t_a^0 \cdot [ 1 + \eta (\frac{\omega}{c_a})^\beta] \; d\omega \nonumber \\
    & = \min \sum_{a \in A_n} \alpha \cdot t_a^0 \cdot [ x_a + \frac{\eta (x_a)^{\beta+1}}{(\beta+1) (c_a)^\beta } ]
\end{align}

At this point, we want to approximate the function $(x_a)^{\beta+1}$ by a piece-wise linear function. For that, we divide the link flow $x_a$ into $K$ segments, between values $\mu_a^{k-1}$ and $\mu_a^{k}$ ; $k \in [1,K]$. We then introduce new decision variables $\lambda_{a,k}^L$ and $\lambda_{a,k}^R$ to activate one of the proportion at a time. The goal is now to find these variables such that :

\[ x_a = \sum_{k=1}^K \mu_a^{k}\lambda_{a,k}^L + \mu_a^{k}\lambda_{a,k}^R\]

Thus, the link travel time function is reformulated as follows : 

\begin{equation}
\label{app1_eq2}
    t_a = t_a^0 \cdot [ 1 + \frac{\eta}{c_a^\beta} \sum_{k=1}^K \mu_a^{k}\lambda_{a,k}^L + \mu_a^{k}\lambda_{a,k}^R ]
\end{equation}

By incorporation equation \ref{app1_eq2} into \ref{app1_eq1}, we transform the program T as the following equivalent program:
\begin{align}
\label{app1_eq3}
    \min \sum_{a \in A_n} & \alpha \cdot t_a^0 \cdot [ x_a + \frac{\eta}{(\beta+1) (c_a)^\beta} \cdot y_a ] \nonumber \\
    \textrm{s.t.}  \nonumber \\
    &y_a =  \sum_{k=1}^K \mu_a^{k}\lambda_{a,k}^L + \mu_a^{k}\lambda_{a,k}^R \quad \forall a \in RN \nonumber \\
    & \sum_{k=1}^K \lambda_{a,k}^L + \lambda_{a,k}^R = 1 \quad \forall a \in RN
\end{align}

\section{User Equilibrium Conditions Equivalency}
\label{app::UE}

In this section, we prove the equivalency of our program to the UE conditions. For that, we formulate the KKT conditions for the continuous setting to show that the solution of our program is at equilibrium. This approach is similar to the work of \cite{du_modeling_2022} and provides the theoretical foundations to formulate such multimodal frameworks and show the applicability of our model to the continuous setting as well. Please note that introducing integer variables does not affect the UE conditions of the program, as stated by \cite{rosenthal_network_1973}, but may lead to differences in the solutions. This is further discussed in the paper.

Let $\mu_{rs} \;( r,s \in E$) denote the dual variable associated with the flow conservation constraint in Equation (\ref{eq7}). Thus, the Lagrangian can be given by: 

\begin{equation}
\label{lagrangian}
L(f, \mu)  = Z + \sum_{r,s \in E} \mu_{r,s} (q^{(r,s)} - \sum_{m \in \Psi} \sum_{p \in P_{rs}^m} f_{p,m}^{(r,s)}  ) 
\end{equation}

We compute the derivatives of the Lagrangian with respect to the path flow to express the KKT conditions. To simplify the equations, we assume in what follows that $\alpha=1$.

1. If $m \notin \{ CP, RS, EH\}$: for all modes, except SMSs, the terms WT, ST, and C in the objective function are constants w.r.t the path flow. Thus, we obtain the following derivatives:

\begin{align}
    \nonumber \\
    \frac{\partial Z}{\partial f_{p,m}^{(i,j)}} & =  \sum_{b \in A} \frac{\partial }{\partial x_b} \{ \sum_{n \in N} \sum_{a \in A_n}  \int_0^{x_a}  t_{a,n}(\omega) \; d\omega \} \times \frac{\partial x_{b}}{\partial f_{p,m}^{(i,j)}} \nonumber \\
    & + \sum_{b \in A} \sum_{m' \in \Psi}  \frac{\partial }{\partial x_{b,m'}} \{ \sum_{a \in A} \sum_{m \in \Psi} [ WT_{m,a} + ST_{m,a} + C_{m,a} ] \cdot x_{a,m} \}  \times \frac{\partial x_{b,m'}}{\partial f_{p,m}^{(i,j)}} \\
    \nonumber \\
    & = \sum_{b \in A} \; t_{b,n}(x_b) \cdot \delta_{b,p,m}^{(i,j)} +  \sum_{b \in A} [ WT_{m,b} + ST_{m,b} + C_{m,b} ] \times \delta_{b,p,m}^{(i,j)}\\
    & = \sum_{b \in A} \; [ t_{b,n}(x_b) + WT_{m,b} + ST_{m,b} + C_{m,b} ] \times \delta_{b,p,m}^{(i,j)}\\
    & = c_{p,m}^{(i,j)}   \qquad \qquad \forall i,j \in E \quad ; \; \forall m \in \Psi - \{ CP, RS, EH\} \quad ; \; \forall p \in P_{ij}^m
\end{align}
\\
2. If $m \in \{ CP, RS, EH\}$: the terms ST and C in the objective function are constants w.r.t the path flow, while WT is not. Let us define $G_b \; (\forall b \in A)$ as below: 
\\
\begin{align}
G_{b,m} &= \frac{\partial WT_{m,b} }{\partial x_{b,m}} \label{eq:28}\\
&= \frac{\partial WT_{m,b}}{\partial f_{p,m}^{(i,j)}} \times \frac{\partial f_{p,m}^{(i,j)}}{\partial x_{b,m}} \label{eq:29}\\ 
& = \frac{1}{R_m} \times \frac{1}{\delta_{b,p,m}^{(i,j)}} \label{eq:30} \\
\nonumber
\end{align}

\noindent The transition from Equation (\ref{eq:28}) to Equation (\ref{eq:29}) is by applying the chain rule of differentiation. The transition from Equation (\ref{eq:29}) to Equation (\ref{eq:30}) is explained by the fact that the waiting time for a passenger using an SMS depends only on the number of passengers using that same service. Also, the incidence matrix  [$(\delta_{a,p,m}^{(i,j)})_{i,j \in E \; ; a \in A \; ; m \in \Psi \; ; p \in P_{ij}^m}$] is predefined and fixed for every network, and depends on neither the flow variables nor the OD demand. Then, as we are working with the same boundaries for functions $f$ and $x$, we can use the property of derivatives, which states that: if $ \frac{\partial x}{\partial f} = g$, then $ \frac{\partial f}{\partial x} = \frac{1}{g}$.
With this, we calculate the derivative of the objective function w.r.t the flow variable.

\begin{align}
    \nonumber \\
    \frac{\partial Z}{\partial f_{p,m}^{(i,j)}} & =  \sum_{b \in A} \frac{\partial }{\partial x_b} \{ \sum_{n \in N} \sum_{a \in A_n}  \int_0^{x_a}  t_{a,n}(\omega) \; d\omega \} \times \frac{\partial x_{b}}{\partial f_{p,m}^{(i,j)}} \nonumber \\
    & + \sum_{b \in A} \sum_{m' \in \Psi}  \frac{\partial }{\partial x_{b,m'}} \{ \sum_{a \in A} \sum_{m \in \Psi} [ WT_{m,a} + ST_{m,a} + C_{m,a} ] \cdot x_{a,m} \}  \times \frac{\partial x_{b,m'}}{\partial f_{p,m}^{(i,j)}} \\
    \nonumber \\
    & = \sum_{b \in A} \; t_{b,n}(x_b) \cdot \delta_{b,p,m}^{(i,j)} + \sum_{b \in A} [ WT_{m,b} + ST_{m,b} + C_{m,b} +  \frac{\partial WT_{m,b} }{\partial x_{b,m}} \cdot x_{b,m} ] \times \delta_{b,p,m}^{(i,j)} \\
    & = \sum_{b \in A} \; [\; t_{b,n}(x_b) + WT_{m,b} + ST_{m,b} + C_{m,b} + G_{b,m} \cdot x_{b,m}\;] \times \delta_{b,p,m}^{(i,j)} \\
    & = c_{p,m}^{(i,j)} \qquad \qquad \forall i,j \in E \quad ; \; \forall m \in \{ CP, RS, EH\} \quad ; \; \forall p \in P_{ij}^m \\
    \nonumber
\end{align}

\noindent In conclusion, the derivative of the Lagrangian is calculated as follows: 

\begin{align}
    \nonumber \\
    \frac{\partial L(f,u)}{\partial f_{p,m}^{(i,j)}} & = \frac{\partial Z}{\partial f_{p,m}^{(i,j)}} + \frac{\partial }{\partial f_{p,m}^{(i,j)}} \{ \sum_{r,s \in E} \mu_{r,s} (q^{(r,s)} - \sum_{m \in \Psi} \sum_{p \in P_{rs}^m} f_{p,m}^{(r,s)} ) \} \\
    & =  c_{p,m}^{(i,j)} - \;\mu_{ij}  \qquad \qquad \forall i,j \in E \quad ; \; \forall m \in \Psi \quad ; \; \forall p \in P_{ij}^m \\
    \nonumber
\end{align}

\noindent Using these derivatives, and following the work of~\cite{sheffi_urban_1984} and~\cite{wang_convex_2021}, we can explicitly formulate the minimization problem in \ref{P1} with the following UE conditions, which complete our proof.

\begin{align}
\nonumber \\
    f_{p,m}^{(i,j)} (c_{p,m}^{(i,j)} - \;\mu_{ij} ) & = 0 \qquad \forall i,j \in E \quad \forall m \in \Psi \quad \forall p \in P_{ij}^m \\
    c_{p,m}^{(i,j)} - \;\mu_{ij}  & \geq 0 \qquad \forall i,j \in E \quad \forall m \in \Psi \quad \forall p \in P_{ij}^m \\
    \sum_{m \in \Psi} \sum_{p \in P_{ij}^m} f_{p,m}^{(i,j)} & = q^{(i,j)} \qquad \forall i,j \in E \\
    f_{p,m}^{(i,j)} & \geq 0 \qquad \forall i,j \in E \quad \forall m \in \Psi \quad \forall p \in P_{ij}^m \\
    \nonumber
\end{align}

\noindent We should note that the generalized cost in case of SMSs, as defined in equation (B.12), includes an additional term that can be seen as a delay induced by the use of SMSs. Specifically, at equilibrium, the commuters choosing a shared mode (carpooling, e-hailing or ridesharing) will experience a delay due to the system-wide definition of waiting time. However, the equivalency conditions state that these commuters will still choose the option with minimum cost, considering this delay, and they cannot change their choices to gain more. Thus, the UE condition remain satisfied.

It is also worth mentioning that if we use a constant value, or a link-based formulation for the waiting time of SMSs, then the additional derivative term in equation (B.12) will not be present and thus, no delay will be experienced, and the solution of our program will satisfy the classical UE principle. However, doing so will not capture the system-wide impact of waiting time on commuters' choices.

Defining The Lagrangian function as in Equation (\ref{lagrangian}) will ensure the correspondence between our objective function (Beckmann Formulation), and the UE conditions in a multi-modal transportation network. Since the variables $f_{p,m}^{(i,j)} \forall i,j \in E $ include all travel options, we only consider the dual variable of the flow conservation constraint. { This means that the solution of our program is at equilibrium, subject to capacity constraints regulating the shared mobility services and public transport. In other words, at the optimal solution, the commuters cannot change their decisions (path and mode choice) to improve their cost while keeping the multi-modal system coherent in terms of capacity and availability of resources. This lead to a definition of equilibrium with side constraints, similar to the problem studied in \cite{larsson_side_1999}. This concept is further discussed in what follows, with an experimental analysis.}

\subsection{\textbf{Equilibrium State' Experimental Discussion}}
Let us consider the network illustrated by \newref{Figure}{UE_net}(A), with one OD pair $(1,2)$. All travel options are available, and the most important parameters are presented in the figure. We relax the PT constraints by considering an unlimited capacity for bus and metro modes. In other words, the PT operates under a fixed schedule and adapts to the passengers' demand. For example, let us consider the bus service. If the buses available are full, bus units (equivalent to three car units) are added to the paths with a high demand. In this test case, the bus mode has the least cost, so all commuters between OD pair $(1,2)$ use the bus to reach their destination. The path flow distribution is shown in \newref{Figure}{UE_net}(B). 90 bus passengers use the first path, while 10 use the second one. The path costs are slightly different, so one could argue that the passengers on the second path will converge to using the first one. However, adding one passenger on the first path will lead to adding one bus unit and, thus, transforms the cost of this path to $14.668$ (higher than the cost of the second path). This means that no passenger can improve the travel costs by changing the path and/or the travel mode corresponding to a UE state.

\begin{figure}[h!]
 \centering
 \includegraphics[width=0.9\textwidth]{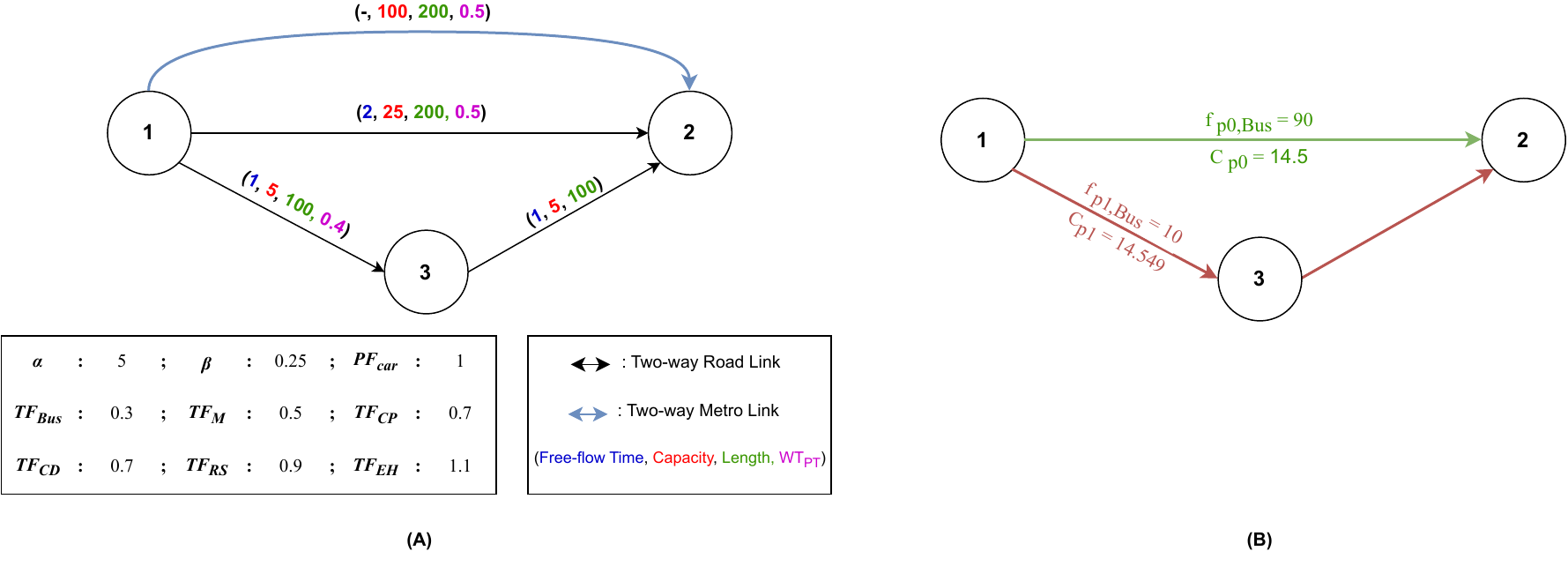}
 \caption{ (A) Toy Network for the Equilibrium State' Experimental Discussion in Multi-Modal Networks. (B) User Equilibrium State between multiple paths of the same travel mode.}
 \label{UE_net}
\end{figure}

To further explore the equilibrium between travel options, we set the various parameters so as to have a similar cost among travel modes, and thus, we focus more on the traffic assignment at equilibrium. \newref{Figure}{UE_net2} shows this assignment. We only represent the paths and modes which are used and have the optimal cost. The bus passengers on the first path have a cost of $15.25$ while commuters using the carpooling service on the second path have a cost of $15.3$. However, if one of these commuters decides to change the path and mode to use the bus service on the first path, a bus unit will be added to the first path, and a car unit will be deleted from the second path. Thus, they will experience a cost of $15.416$ while the carpooling generalized cost on the second path will then be $14.8$. This shows that, for this scenario, the solution returned by our model is at user equilibrium.

\begin{figure}[h!]
 \centering
 \includegraphics[width=0.9\textwidth]{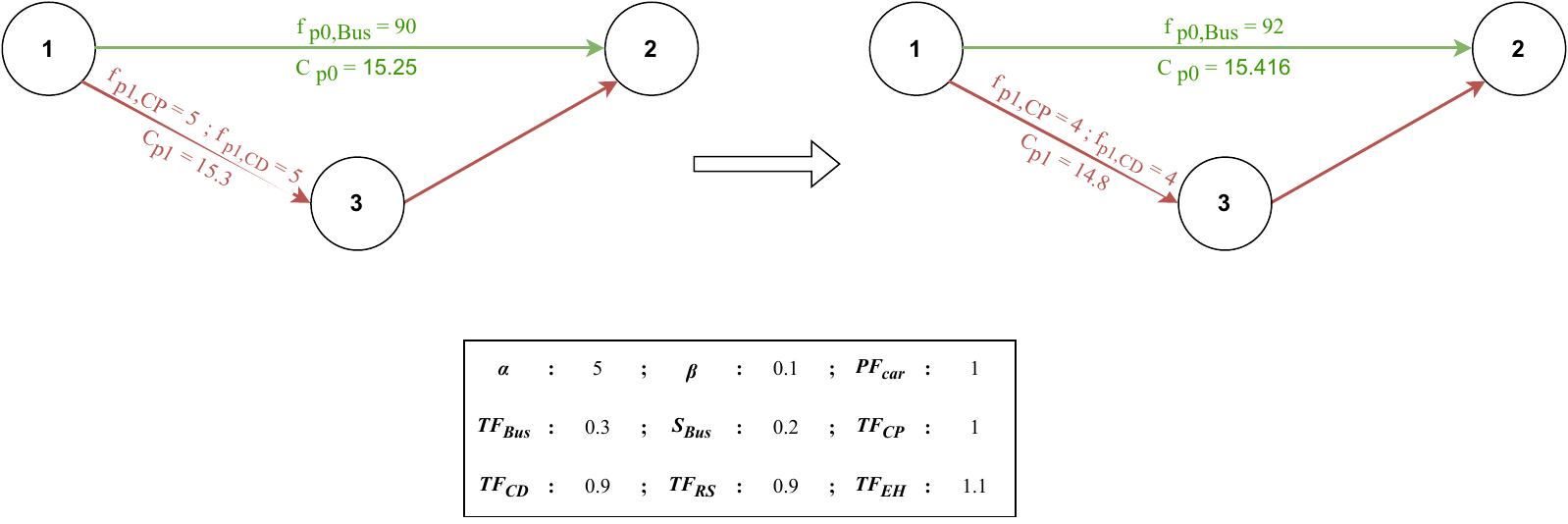}
 \caption{User Equilibrium State between multiple travel modes.}
 \label{UE_net2}
\end{figure}

\section{Upper Bound for the Price of Anarchy}
\label{app::UB_PoA}
Let us consider the PoA formula described by Equation (\ref{Eq: PoA}). The goal is to bind transportation inefficiencies with selfish routing by a maximum value representing the worst-case scenario. Since the PT adheres to a pre-defined schedule with a fixed configuration and operates independently of the passenger flow, and since we do not consider any congestion effect in the metro network, nor for the walking and biking networks, the following part is only valid for the road network, with BPR as a travel time function. Consequently, we omit the network index in the link travel time function in what follows.

{In the work of ~\cite{qiao_price_2023}, the authors provided an upper bound for the BPR function, exclusively considering the road network. In our case, we consider multimodal networks, with travel costs including other terms in addition to the travel time. Thus, this upper bound cannot be directly applied to our study. However, we will use it as a starting point in what follows. }

It is proven that ($\beta$ being the power in the BPR function): 
\begin{equation}
    \sum_{a \in A} t_a(x_a^{UE}) \leq  \frac{(\beta+1)^{1+\frac{1}{\beta}}}{(\beta+1)^{1+\frac{1}{\beta}} - \beta} \cdot \sum_{a \in A} t_a(x_a^{SO})
\end{equation}
Let us represent this upper bound for the BPR PoA with $B(\beta)$. We then multiply both sides of the inequality by the value of time.
\begin{equation} \label{Inq:2}
    \sum_{a \in A} \alpha \cdot t_a(x_a^{UE}) \leq  B(\beta) \cdot \sum_{a \in A} \alpha \cdot t_a(x_a^{SO}) 
\end{equation}
Additionally, let $F(x,q)$ be defined as follows: 
\begin{align}
    F(x,q) & = \sum_{a \in A} \sum_{m \in \Psi} 
    ( \alpha \frac{q_m}{R_m} \cdot b_m + \alpha ST_{m,a} + C_{m,a}) \cdot x_{a,m} \nonumber \\
    & = \sum_{a \in A} \sum_{m \in \Psi}  \alpha \frac{q_m}{R_m} \cdot b_m \cdot x_{a,m} + \sum_{a \in A} \sum_{m \in \Psi} H_{m,a} \cdot x_{a,m}
\end{align}
with: 
\begin{align}
     & b_m = 
     \begin{cases}
           1  & if  \; m \in \{CP,EH,RS\} \\
           0 & otherwise \\ 
     \end{cases}  \nonumber \\
     & H_{m,a} = \alpha ST_{m,a} + C_{m,a}  \nonumber
\end{align}

\noindent We add $F(x^{UE}, q^{UE}) + B(\beta) \cdot F(x^{SO}, q^{SO})$ in both sides of Inequality (\ref{Inq:2}). We now have:
\begin{equation} \label{Inq:3}
    \sum_{a \in A} \alpha \cdot t_a(x_a^{UE}) + F(x^{UE}, q^{UE}) + B(\beta) \cdot F(x^{SO}, q^{SO})  \leq  B(\beta) \cdot \sum_{a \in A} \alpha \cdot t_a(x_a^{SO}) + F(x^{UE}, q^{UE}) + B(\beta) \cdot F(x^{SO}, q^{SO}) 
\end{equation}

\noindent From Inequality (\ref{Inq:3}), we can deduce: 
\begin{align} \label{Inq:4}
    \sum_{a \in A} \alpha \cdot t_a(x_a^{UE}) + F(x^{UE}, q^{UE}) & \leq  B(\beta) \cdot [ \sum_{a \in A} \alpha \cdot t_a(x_a^{SO}) + F(x^{SO}, q^{SO}) ] + F(x^{UE}, q^{UE}) \nonumber \\ 
    \frac{\sum_{a \in A} \alpha \cdot t_a(x_a^{UE}) + F(x^{UE}, q^{UE}) }{\sum_{a \in A} \alpha \cdot t_a(x_a^{SO}) + F(x^{SO}, q^{SO})} & \leq  B(\beta) + \frac{F(x^{UE}, q^{UE})}{\sum_{a \in A} \alpha \cdot t_a(x_a^{SO}) + F(x^{SO}, q^{SO})}
\end{align}

\noindent Now, we aim to provide an upper bound for the right side of Inequality (\ref{Inq:4}). 

\begin{equation}
    \frac{F(x^{UE}, q^{UE})}{\sum_{a \in A} \alpha \cdot t_a(x_a^{SO}) + F(x^{SO}, q^{SO})} = \frac { \sum_{a \in A} \sum_{m \in \Psi}  \alpha \frac{q_m^{UE}}{R_m} \cdot b_m \cdot x_{a,m}^{UE} + \sum_{a \in A} \sum_{m \in \Psi} H_{m,a} \cdot x_{a,m}^{UE}  }  {\sum_{a \in A} \alpha \cdot t_a(x_a^{SO}) + \sum_{a \in A} \sum_{m \in \Psi}  \alpha \frac{q_m^{SO}}{R_m} \cdot b_m \cdot x_{a,m}^{SO} + \sum_{a \in A} \sum_{m \in \Psi} H_{m,a} \cdot x_{a,m}^{SO}} \nonumber
\end{equation}

\noindent Here, we assume that $\sum_{a \in A} \sum_{m \in \Psi} H_{m,a} \geq 1 $. In other words, for every link, the sum of the service times and monetary costs for all the available modes are assumed to be higher than one.

\begin{equation}
    \frac{F(x^{UE}, q^{UE})}{\sum_{a \in A} \alpha \cdot t_a(x_a^{SO}) + F(x^{SO}, q^{SO})} \leq \sum_{a \in A} \sum_{m \in \Psi}  \alpha \frac{q_m^{UE}}{R_m} \cdot b_m \cdot x_{a,m}^{UE} + \sum_{a \in A} \sum_{m \in \Psi} H_{m,a} \cdot x_{a,m}^{UE}
\end{equation}

\noindent From the flow conservation condition defined by Constraint (\ref{eq7}) in the proposed model, we know that, for every feasible solution, and the UE solution in particular, the sum of all commuters in the network equals the total demand. Thus, the number of SMS passengers can equal, at most, the total demand. In other words:

\begin{align}
    \sum_{m \in \Psi} q_m^{UE} \leq Q \nonumber \\
    \sum_{m \in \Psi} \frac{q_m^{UE}}{R_m} \leq \frac{Q}{R_m} \leq Q 
\end{align}

\noindent Additionally, since the path flow of every mode cannot exceed the total demand, the link flow is also bounded by $Q$.

\begin{align} \label{Inq:5}
    &x_{a,m} \leq Q & \forall a \in A \nonumber \\ 
    &\sum_{m \in \Psi} \frac{q_m^{UE}}{R_m} \cdot x_{a,m}  \leq Q^2 & \forall a \in A \nonumber \\ 
    &\sum_{a \in A} \sum_{m \in \Psi} \frac{q_m^{UE}}{R_m} \cdot x_{a,m} + \sum_{a \in A} \sum_{m \in \Psi} H_{m,a} \cdot x_{a,m}^{UE} \leq \sum_{a \in A} Q^2 + \sum_{a \in A} \sum_{m \in \Psi} H_{m,a} \cdot Q \nonumber \\
    &\sum_{a \in A} \sum_{m \in \Psi} \frac{q_m^{UE}}{R_m} \cdot x_{a,m} + \sum_{a \in A} \sum_{m \in \Psi} H_{m,a} \cdot x_{a,m}^{UE} \leq |A| \cdot Q^2 + |A||\Psi | \cdot Q 
\end{align}

\noindent From Inequalities (\ref{Inq:4}) and (\ref{Inq:5}), we have: 

\begin{equation}
    \frac{\sum_{a \in A} \alpha \cdot t_a(x_a^{UE}) + F(x^{UE}, q^{UE}) }{\sum_{a \in A} \alpha \cdot t_a(x_a^{SO}) + F(x^{SO}, q^{SO})} \leq  B(\beta) + |A| \cdot Q^2 + |A||\Psi | \cdot Q 
\end{equation}



\newpage


\setcounter{table}{0}
\setcounter{figure}{0}

\end{document}